\documentclass[a4paper]{amsart}

\usepackage{amssymb,pstricks,diagrams, amscd}
\usepackage[totalwidth=17cm,totalheight=24cm]{geometry}
\usepackage{graphicx}
\input psfig
\diagramstyle[nohug]

\begin{document}

\newtheorem{cor}{Corollary}[section]
\newtheorem{theorem}[cor]{Theorem}
\newtheorem{prop}[cor]{Proposition}
\newtheorem{lemma}[cor]{Lemma}
\theoremstyle{definition}
\newtheorem{defi}[cor]{Definition}
\theoremstyle{remark}
\newtheorem{remark}[cor]{Remark}
\newtheorem{example}[cor]{Example}

\newcommand{\cD}{{\mathcal D}}
\newcommand{\cM}{{\mathcal M}}
\newcommand{\cT}{{\mathcal T}}
\newcommand{\cML}{{\mathcal M\mathcal L}}
\newcommand{\cGH}{{\mathcal G\mathcal H}}
\newcommand{\C}{{\mathbb C}}
\newcommand{\N}{{\mathbb N}}
\newcommand{\R}{{\mathbb R}}
\newcommand{\Z}{{\mathbb Z}}
\newcommand{\Kt}{\tilde{K}}
\newcommand{\Mt}{\tilde{M}}
\newcommand{\dr}{{\partial}}
\newcommand{\kappab}{\overline{\kappa}}
\newcommand{\pib}{\overline{\pi}}
\newcommand{\Sigmab}{\overline{\Sigma}}
\newcommand{\gd}{\dot{g}}
\newcommand{\diff}{\mbox{Diff}}
\newcommand{\dev}{\mbox{dev}}
\newcommand{\devb}{\overline{\mbox{dev}}}
\newcommand{\devt}{\tilde{\mbox{dev}}}
\newcommand{\vol}{\mbox{Vol}}
\newcommand{\hess}{\mbox{Hess}}
\newcommand{\db}{\overline{\partial}}
\newcommand{\Sigmat}{\tilde{\Sigma}}

\newcommand{\cunc}{{\mathcal C}^\infty_c}
\newcommand{\cun}{{\mathcal C}^\infty}
\newcommand{\dd}{d_D}
\newcommand{\dmin}{d_{\mathrm{min}}}
\newcommand{\dmax}{d_{\mathrm{max}}}
\newcommand{\Dom}{\mathrm{Dom}}
\newcommand{\dn}{d_\nabla}
\newcommand{\ded}{\delta_D}
\newcommand{\delmin}{\delta_{\mathrm{min}}}
\newcommand{\delmax}{\delta_{\mathrm{max}}}
\newcommand{\hmin}{H_{\mathrm{min}}}
\newcommand{\maxi}{\mathrm{max}}
\newcommand{\oL}{\overline{L}}
\newcommand{\oP}{{\overline{P}}}
\newcommand{\Ran}{\mathrm{Ran}}
\newcommand{\tgamma}{\tilde{\gamma}}
\newcommand{\cotan}{\mbox{cotan}}
\newcommand{\lambdat}{\tilde\lambda}
\newcommand{\St}{\tilde S}

\newcommand{\II}{I\hspace{-0.1cm}I}
\newcommand{\III}{I\hspace{-0.1cm}I\hspace{-0.1cm}I}
\newcommand{\note}[1]{\marginpar{\tiny #1}}

\title{AdS manifolds with particles and earthquakes on singular surfaces}
\author{Francesco Bonsante}
\address{Scuola Normale Superiore\\
Piazza dei Cavalieri\\
56100 Pisa, Italy}
\email{bonsante@sns.it}
\author{Jean-Marc Schlenker}
\address{Institut de Math\'ematiques, UMR CNRS 5219\\
Universit\'e Toulouse III\\
31062 Toulouse cedex 9, France}
\email{schlenker@math.ups-tlse.fr}
\date{June 2007 (v3)}

\begin{abstract}
We prove two related results. The first is an ``Earthquake Theorem'' 
for closed hyperbolic surfaces with cone
singularities where the total angle is less than $\pi$: any two such
metrics in are connected by a unique left earthquake. 
The second result is that the space of
``globally hyperbolic'' AdS manifolds with ``particles'' -- 
cone singularities (of given angle) along time-like
lines -- is parametrized by the product of two copies of the Teichm\"uller
space with some marked points (corresponding to the cone singularities).
The two statements are proved together.
\end{abstract}
\maketitle


\section{Introduction and results}

\subsection*{The Earthquake Theorem.}

Let $\Sigma$ be a closed surface with a hyperbolic metric $h$, and let
$\lambda$ be a measured lamination on $\Sigma$. Then $\lambda$ can be uniquely
realized as a measured geodesic lamination for $h$. Thurston
\cite{thurston-notes,thurston-earthquakes} 
defined the image of $h$ by the {\it right earthquake}
along $\lambda$, called $E^r_\lambda(h)$ here, in a way which can be described
simply when the support of 
$\lambda$ is a disjoint union of closed curves: it is obtained by cutting
$\Sigma$ along each geodesics in the support of $\lambda$, doing a 
fractional Dehn twist
by the length corresponding to the weight assigned to the curve by $\lambda$,
and gluing back. This defines a map:
$$ E^r:\cML_\Sigma\times \cT_\Sigma\rightarrow \cT_\Sigma~. $$

Thurston then proved that the corresponding action of $\cML_\Sigma$ on
$\cT_\Sigma$ is simply transitive: given $h, h'\in \cT_\Sigma$, there is a
unique $\lambda\in \cML_\Sigma$ such that $E^r_\lambda(h)=h'$. A different
proof was given by Kerckhoff \cite{kerckhoff}.

\subsection*{The Mess proof of the Earthquake Theorem.}

Yet another proof of the Earthquake Theorem was later discovered by 
Mess \cite{mess} as a by-product of the geometric properties of globally
hyperbolic maximal compact (GHMC) Anti-de Sitter (AdS) manifolds. Mess
discovered that such
manifolds share several remarkable properties with quasi-Fuchsian hyperbolic
manifolds. In particular those manifolds (containing a space-like
surface diffeomorphic to $\Sigma$) are uniquely determined by two 
hyperbolic metrics on $\Sigma$, called their ``left'' and ``right''
representations, which are analogs in the AdS context of the conformal
metrics at infinity for quasifuchsian manifolds. This result of Mess
can be interpreted as an analog of the Bers double uniformization theorem
for quasifuchsian 3-manifolds.

GHMC AdS manifolds also have a convex core, with a boundary 
which has a hyperbolic induced metric and which is ``pleated'' along a 
measured geodesic laminations. Earthquakes are natural analogs in AdS geometry
of grafting in quasifuchsian geometry, and some geometric properties
of GHMC AdS manifold then yield the Earthquake Theorem as a consequence
of the Mess parametrization of the space of those manifolds by
$\cT_\Sigma\times \cT_\Sigma$.

\subsection*{Surfaces with cone singularities.}

From this point on we consider a closed surface $\Sigma$ with 
$n$ distinct marked points $x_1, \cdots, x_n$. We are interested
in the hyperbolic metrics on $\Sigma$ with cone singularities
at the $x_i$. Given such a metric, we call $\theta_i$ the angle at $x_i$,
$1\leq i\leq n$. It follows from a result of Troyanov \cite{troyanov} that,
given the $\theta_i$, those metrics are in one-to-one correspondence with the
conformal structures on $\Sigma$, so that, considered up to the isotopies
fixing the $x_i$, those metrics are parametrized by $\cT_{\Sigma,n}$, the
Teichm\"uller space of $\Sigma$ with $n$ marked points. Setting
$\theta:=(\theta_1, \cdots, \theta_n)$, we will often denote by
$\cT_{\Sigma, n,\theta}$ the space of hyperbolic metrics on $\Sigma$ with cone
singularities at the $x_i$ of angles given by the $\theta_i$ (considered up to
isotopies fixing the marked points). 

It is interesting to note at this point that the theory of geodesic
laminations works on hyperbolic surfaces with cone singularities quite like it
does on closed hyperbolic surfaces, as long as the cone angles $\theta_i$ are
less than $\pi$.

\begin{lemma} \label{lm:realization}
Suppose that $\theta_i<\pi, 1\leq i\leq n$. Then each measured lamination in
the complement of the $x_i$ in $\Sigma$ can be realized uniquely as a geodesic
lamination. 
\end{lemma}

The proof (which is elementary) can be found in section 3. In addition, still
under the hypothesis that the angles are less than $\pi$, the geodesic
laminations can not come too close to the singular points, and it follows
that earthquakes along measured geodesic laminations do not change the angles
at the cone points.

We will call $\cML_{\Sigma, n}$ the space of measured laminations on the
complement of the $x_i$ in $\Sigma$. It is well known (see \cite{FLP}) 
that $\cML_{\Sigma, n}$ is homeomorp¶ic to a ball 
of dimension $6g-6+2n$, where $g$ is
the genus of $\Sigma$.

\subsection*{Earthquakes on singular surfaces.}

The first result of this paper is an extension of the Earthquake Theorem to
hyperbolic surfaces with cone singularities.

\begin{theorem} \label{tm:earthquake}
For all $h,h'\in \cT_{\Sigma,n,\theta}$, there is a
unique $\lambda\in \cML_{\Sigma,n}$ such that $E_r(\lambda)(h)=h'$. 
\end{theorem}

\subsection*{AdS manifolds with particles.}

The second theme considered here --- which is strongly related to the first
--- concerns 3-dimensional AdS manifolds with cone singularities along
time-like lines. Such cone singularities are called ``particles'' here, since
they are sometimes used in the $2+1$-gravity 
community to model massive, spin-less point particles. A precise definition is
given in section 3. We are in particular interested in ``globally hyperbolic
compact maximal'' AdS manifolds with ``particles'', extending those considered
by Mess \cite{mess}. An AdS manifold with ``particles'' (cone
singularities along time-like lines) $M$ is GHMC if:
\begin{itemize}
\item it contain a closed, oriented, locally convex
  space-like surface $S$ which is ``orthogonal to the
  singular line'' (in a manner which is described in section 3),
\item every inextendible time-like curve in $M$ intersects $S$ exactly once,
\item if $M'$ is another AdS manifold with particles satisfying the first two
  properties in which $M$ can be isometrically embedded, then $M'=M$.
\end{itemize}

\begin{defi}
We call $\cGH_{\Sigma,n}$ the space of GHMC
AdS metrics on $\Sigma\times \R$, with cone singularities at the lines
$\{ x_i\}\times \R$ for $1\leq i\leq n$, considered up to isotopies fixing
the singular lines. Given $\theta:=(\theta_1, \cdots, \theta_n)\in (0,\pi)^n$,
we also call $\cGH_{\Sigma,n,\theta}$ the subspace of those metrics for which
the angle at the line $\{ x_i\}\times \R$ is $\theta_i$.
\end{defi}

Note that we will sometimes abuse notations and write about a GHMC AdS metric
or a GHMC AdS manifold indifferently.

\subsection*{The right and left metrics associated to a GHMC manifold.}

Let $\theta:=(\theta_1, \cdots, \theta_n)\in (0,\pi)^n$, and let 
$M\in \cGH_{\Sigma, n,\theta}$. By definition $M$ contains an oriented,  
space-like, convex
surface $S$ which is ``orthogonal'' to the singular lines. Let $I$ be the
induced metric on $S$, let $J$ be the complex structure associated to $I$, and
let $B$ be the shape operator of $S$. It is then possible to define two
metrics $\mu_l$ and $\mu_r$ on $S$ as:
$$ \mu_l := I((E+JB)\cdot, (E+JB)\cdot)~, ~~ \mu_r := I((E-JB)\cdot,
(E-JB)\cdot)~. $$
This corresponds to the metrics $I^\#_\pm$ defined in \cite{minsurf}, the
notations presented here are better suited for our needs. It is proved in
\cite{minsurf} that those two metrics are hyperbolic, with cone singularities
at the intersections of $S$ with the ``particles'', where their angle is equal
to the angle of $M$ at the corresponding ``particle''. Moreover it is also
proved in \cite{minsurf} that $\mu_l$ and $\mu_r$ (considered up to isotopy)
do not depend on the choice of $S$. So this construction defines two maps
$$ \mu_l, \mu_r : \cGH_{\Sigma,n,\theta}\rightarrow \cT_{\Sigma,n,\theta}~. $$
When no ``particle'' is present, $\mu_l$ and $\mu_r$ are the hyperbolic
metrics corresponding to the right and left representations considered by Mess
\cite{mess}. 

The second main result here is that, as for non-singular GHMC AdS manifolds,
the maps $\mu_l$ and $\mu_r$ provide a parametrization of
$\cGH_{\Sigma,n,\theta}$. 

\begin{theorem} \label{tm:mess}
The map $(\mu_l, \mu_r): 
\cGH_{\Sigma,n,\theta}\rightarrow \cT_{\Sigma,n,\theta}\times
\cT_{\Sigma,n,\theta}$ is one-to-one.
\end{theorem}

This statement can be construed as an extensions of Mess' AdS version
of the Bers double uniformization theorem to AdS manifolds with ``particles''.
Note that on the hyperbolic side an extension of the Bers double uniformization
theorem to quasifuchsian manifolds with ``particles'' -- cone singularities
of angle less than $\pi$ along lines going from one boundary at infinity to
the other -- might well hold but it has not been proved yet (a first
step is made in \cite{qfmp}).

\subsection*{The structure of GHMC AdS manifolds with particles.}

The proof of Theorem \ref{tm:earthquake} uses Theorem \ref{tm:mess}
along with some properties which were
discovered by Mess for non-singular GHMC AdS manifolds, which extend
directly to GHMC AdS manifolds with particles. In particular, 
those manifolds contain a smallest convex subset --- where ``smallest''
is understood with respect to the inclusion --- called its convex 
core. We call $CC(M)$ the convex core of a GHMC AdS manifold $M$.

\begin{lemma} \label{lm:cc}
The boundary of $CC(M)$ has two connected components, called 
$\dr_+ CC(M)$ and $\dr_- CC(M)$. Each is a space-like surface (outside its
intersection with the singular set of $M$), which is ``orthogonal'' to $M_s$
(as defined in section 3). Its induced metric is hyperbolic, and it is
``pleated'' along a measured geodesic lamination.
\end{lemma}

The structure of $M$ can be readily understood from its convex core, exactly
as for non-singular GHMC AdS manifolds.

\begin{lemma} \label{lm:structure}
  \begin{enumerate}
  \item Let $x\in \dr_- CC(M)$, and let $H$ be a space-like plane containing 
$x$ which is a support plane of $CC(M)$ at $x$. Let $n$ be the
future-oriented unit vector at $x$ which is orthogonal to $H$.
Then the geodesic maximal segment starting from $x$ in the direction of $n$
has length $\pi/2$. 
\item 
$M$ is the union of the future of $\dr_-CC(M)$ and of the
past of $\dr_+CC(M)$, and their intersection is equal to $CC(M)$.
  \end{enumerate}
\end{lemma}

Moreover, the metric on the future of the past boundary of the convex core can
expressed, in a fairly simple way, in terms of the induced metric and the
measured bending lamination, on the past boundary of the convex core.

\subsection*{Outline of the proofs}

The basic idea of the proof of Theorem \ref{tm:earthquake} is to use a ``deformation''
argument, in which we fix $\theta\in (0,\pi)^n$. 
We fix a hyperbolic $h\in \cT_{\Sigma,n,\theta}$ and consider the map 
$$ 
\begin{array}{rccc}
E^r_\cdot (h) : & \cML_{\Sigma,n} & \rightarrow & \cT_{\Sigma,n,\theta} \\
& \lambda & \mapsto & E^r_\lambda(h)~. $$
\end{array}
$$
Our goal is to show that this map is homeomorphism. This follows from some basic
points:
\begin{enumerate}
\item $E^r_\cdot (h)$ is a local homeomorphism,
\item it is proper,
\item its target space $\cT_{\Sigma,n,\theta}$ is simply connected, while $\cML_{\Sigma,n}$
is connected.
\end{enumerate}

A key point of the proof, however, is to use the settings of both 
Theorem \ref{tm:earthquake} and Theorem \ref{tm:mess}. The first point
is to prove the equivalence between those two statements; the fact
that $E^r_\cdot (h)$ is a local homeomorphism can be proved in the setting of Theorem 
\ref{tm:mess}, while the fact that $E^r_\cdot (h)$ is proper can
be shown fairly easily on the side of Theorem \ref{tm:earthquake}, where it
appears as the following compactness statement.

\begin{lemma} \label{lm:compact}
Let $\theta=(\theta_1,\cdots, \theta_n)$ be fixed. Let $\mu\in
\cT_{\Sigma,n,\theta}$. The map $E^r_\cdot(\mu):\cML_{\Sigma,n}\rightarrow
\cT_{\Sigma,n,\theta}$ is proper.
\end{lemma}

The equivalence between Theorem \ref{tm:earthquake} and Theorem
\ref{tm:mess} uses crucially the existence of a ``convex core'' in
the GHMC AdS manifolds with particles considered here, as proved 
in section 5. It is then possible to prove, in section 6, that
the two main theorems are equivalent, using the relations between
the representation of the fundamental group of a GHMC AdS manifold
with particles, on one hand, and 
the induced metrics and measured lamination on the boundary of the
convex core, on the other.

\subsection*{Acknowledgements.} We would like to thank Thierry Barbot for
several useful remarks and comments.

\section{Background material}\label{sec:back}

This section contains a number of important facts which are presumably
well-known (or close to facts which are well-known), mostly concerning
hyperbolic surfaces or AdS manifolds.

\subsection*{Geodesic laminations on cone surfaces.}

It was mentioned in the introduction that the theory of geodesic laminations
work well on hyperbolic surfaces with cone singularities, provided that the
singular angles are less than $\pi$. A key reason for this is that embedded
geodesics can not come too close to the singularities.

\begin{remark}
There is a decreasing function $\rho:(0,\pi)\rightarrow \R_{>0}$ as follows. Let
$S$ be a surface with a complete hyperbolic metric with cone singularities,
with positive singular curvature,
and let $x$ be one of the singular points, with cone angle $\theta\in
(0,\pi)$. Then any complete embedded geodesic in $S$ is at distance at least
$\rho(\theta)$ from $x$. 
\end{remark}

The proof is left to the reader, the basic idea is that, in a neighborhood of
$x$, the metric on $S$ is isometric to the metric on a hyperbolic cone of
angle $\theta$ near the vertex. However, in such cones, the geodesics which
are too close from the vertex can not be embedded.

\begin{lemma}
Let $S$ be a closed hyperbolic surface with cone singularities, with angle in
$(0,\pi)$ at the singular points. Let $S_0$ be the complement of the singular
points. Then any lamination in $S_0$ can be realized as a geodesic lamination
for the hyperbolic metric on $S_0$.  
\end{lemma}

\begin{proof}
Let $x_1, \cdots, x_n$ be the singular points on $S$, and let $g$ be the
hyperbolic metric on the complement of the $x_i$. 
It follows for instance from \cite{troyanov} that there exists a one-parameter
family of hyperbolic metrics $(g_t)_{t\in [0,1]}$ with cone singularities (or
cusps) at the $x_i$ such that $g_0=g$ and that $g_1$ is a complete hyperbolic
metric with cusps at the $x_i$. Moreover it is possible to demand that the
angle at the cone singularities remain in $(0,\pi-\epsilon)$ for all $t\in
(0,1)$, where $\epsilon$ is some strictly positive constant.

Let $\lambda$ be a lamination on $S_0$. It is well-known (see \cite{FLP}) that
it can be realized uniquely as a lamination which is geodesic for $g_1$. 
Let $E$ be the set of all $t\in [0,1]$ such that, for all $s\in [t,1]$, 
$\lambda$ can be realized uniquely as a geodesic lamination $\lambda_s$ in
$g_s$. Then $E$ is open (because geodesic laminations can be deformed
to ``follow'' a deformation of the underlying hyperbolic metric) and
closed by the previous remark. So $0\in E$, which proves the statement.
\end{proof}

\subsection*{Particles in AdS manifolds.}

The Anti-de Sitter space $H^3_1$ --- called the AdS space here --- is a Lorentz
space of constant curvature $-1$. We only consider it here in dimension 3. 
It is not simply connected. It can be obtained as the quadric:
$$ H^3_1 := \{ x\in \R^4_2 ~ |  ~ \langle x,x\rangle =-1\} $$
in the $4$-dimension flat space of signature $2$, $\R^4_2$, with the induced
metric. 

There is a useful projective model of the AdS space, obtained by considering
the interior of a quadric of signature $(1,1)$ in the sphere $S^3$, with its
``Hilbert metric''. We do not elaborate much here on the geometry of the 
AdS space and refer the reader to e.g. \cite{O,mess}.

\begin{remark}
Let us just remark that there is a natural identification of the isometry group of 
$H^3_1$ with $PSL(2,\mathbb R)\times PSL(2,\mathbb R)$. Moreover there is an isometric 
embedding of the hyperbolic plane $\mathbb H^2$ in $H^3_1$ such that

1) the image is a spacelike geodesic plane.

2) the embedding is $PSL(2,\mathbb R)$ equivariant, where $PSL(2,\mathbb R)$ acts on $H^3_1$ by the diagonal
action.

We call such embedding the \emph{standard} embedding.
\end{remark}

An AdS manifold is a manifold locally modeled on the AdS space. We are
interested here in AdS cone-manifolds of a special kind, which have 
``conical'' singularities along time-like lines; such singularities will
be called ``particles''. They have a simple local description. Consider
a time-like geodesic $c$ in the AdS space, and let $D$ be a domain
in $H^3_1$ bounded by two geodesic time-like half-planes, both bounded by $c$,
with an angle $\theta$ between them (the angle is well-defined since $c$ is
supposed to be time-like). There is a natural way to glue isometrically
the two half-planes so that the identification is the identity on $c$,
and we call $H^3_{1,\theta}$ the resulting space.
The complement of $c$ in $H^3_{1,\theta}$ is locally modeled
on the AdS space --- there is no singularity at the gluing --- while $c$
corresponds to a cone singularity with angle $\theta$, which is the local
model we use for ``particles''.

We then define an {\it AdS manifold with particles} to be a manifold such that
the complement of a disjoint union of open curves is endowed with a Lorentz
metric, and such that each point has a neighborhood isometric either to an
open subset of the AdS space or to a neighborhood of a point of $c$ in 
$H^3_{1,\theta}$, for some value of $\theta$ in $(0,\pi)$.

Note that by construction the angle $\theta$ is constant along a ``particle''. 

\subsection*{GHMC AdS manifolds with particles.}

In the local model described above for the neighborhood of the ``particles'' 
there is a natural notion of ``horizontal'' totally geodesic plane: those are
the image, under the gluing construction, of the restriction to $D$ of the
totally geodesic planes orthogonal to $c$ in $H^3_1$. Note that, under the
gluing construction, thoses planes are sent to totally geodesic surfaces,
i.e., the intersection of those planes with the two half-planes bounding $D$
are glued together and no singularity occurs (except the cone point). By construction there is a
unique horizontal plane containing each point of $c$.

We define a {\it space-like surface} in $H^3_{1,\theta}$ to be a surface
which: 
\begin{itemize}
\item intersects the singular set $c$ at exactly one point $x$,
\item is space-like outside $x$ and locally convex,
\item is such that the tangent plane at a sequence of points converging to $x$
  converges to the tangent plane to the (unique) horizontal plane at $x$.
\end{itemize}
In other terms, what we call ``space-like surfaces'' are really surfaces
which, in addition to being ``space-like'' outside the singular set, are
``orthogonal'' to the singular set in a natural way. 

Returning to an AdS manifold with particles $M$, recall that any point in $M$
has a neighborhood which is isometric to an open set in some $H^3_{1,\theta}$,
for some $\theta\in (0,\pi)$. It is therefore quite natural to define a {\it
  space-like} surface in $M$ as a subset which is corresponds, in each of the
neighborhoods, to a space-like surface in $H^3_{1,\theta}$.

There is also a natural notion of {\it time-like curves} in AdS manifolds with
particles; they are curves which are time-like (in the usual sense) outside
the singular set, and which might follow segments of the singular set, but in
such a way that any time function is monotonous along them.

\subsection*{Globally hyperbolic AdS manifolds.}

We are now almost ready to define a GHMC AdS manifold with particles. Given an
AdS manifold with particles $M$, a {\it Cauchy surface} in $M$ is a closed,
space-like surface such that any inextendible time-like curve in $M$
intersects $S$ exactly once. 

\begin{prop}
If $\Sigma$ is a Cauchy surface of a AdS manifold with particles $M$, then
topologically $M=\Sigma\times\mathbb R$.
\end{prop}

This proposition can be proved by adapting the general result for globally
hyperbolic spacetimes without singularities to the case with particles: one
construct a timelike vector-field on $M$ that is tangent to the singular
locus. The flow of such a field restricted to $S$ realizes a diffeomorphism
between a regular neighbourhood of $S\times\{0\}$ in $S\times\mathbb R$ and
$M$.

\begin{defi} \label{df:convex}
Let $M$ be an AdS manifolds with particles. It is {\bf convex 
globally hyperbolic}
(called GH here) if it contains a Cauchy surface which is locally 
convex. It is {\bf convex globally hyperbolic 
  maximal} (or GHM) if moreover any GH AdS manifold with particles $M'$
containing a subset isometric to $M$ is itself isometric to $M$.
\end{defi}

\begin{prop}\label{ext:prop}
Let $M$ be an AdS GH spacetime with particles. 
There exists a GHM AdS spacetime with particles $M'$, called the maximal
extension of $M$ , in which
$M$ isometrically embeds. Moreover, two maximal extensions are isometric.
\end{prop}

\begin{proof}[Sketch of the proof]
The proof of this proposition follows the same steps of the analog result for
GH spacetimes without singularities~\cite{Choquet-Bruhat}. Actually the way to
adapt the original proof of~\cite{Choquet-Bruhat} to the case without
singularity was suggested to the authors by Thierry Barbot.

The existence of the maximal extension is an application of Zorn's Lemma. We
consider on the set of extensions of $M$ the order given by the isometric
inclusion.  Such an order turns to be inductive so a maximal element exists.

The uniqueness of the maximal extension is more delicate. One proves that given
two extensions $M_1, M_2$, there is another extensions containing both of
them. The idea is to consider the pairs $(N_1,N_2)$ such that
\begin{enumerate}
\item $N_i$ is a GH spacetime contained in $M_i$ and containing $M$.
\item the isometry of $M\rightarrow M_2$ extends to an isometry
  $N_1\rightarrow M_2$ sending $N_1$ to $N_2$.
\end{enumerate}

Clearly there is a natural order on such pairs.  Again an application of
Zorn's Lemma ensures that there is a maximal element among those pairs, say
$(N_1,N_2)$. The idea is to consider the space $\hat M$ obtained by gluing
$M_1$ and $M_2$ identifying $N_1$ to $N_2$. If we prove that $\hat M$ is
topologically a manifold, then it is clear that the Lorentzian structures of
$M_1$ and $M_2$ induce a Lorentzian structure on $\hat M$ in such a way that
the projection $M_1\cup M_2\rightarrow\hat M$ is an isometric embedding on
each component. Moreover simple arguments show that a Cauchy surface of $M$ is
a Cauchy surface for $\hat M$, thus $\hat M$ is GH extension of $M$ containing
both $M_1$ and $M_2$.

Let us prove that $\hat M$ is a manifold. Since the projection $\pi$ is open,
every point has a neighbourhood homeomorphic to $\mathbb R^3$. The only point
to check is that $\hat M$ is a Hausdorff space.  By contradiction, assume
there exist points $x,y\in\hat M$ whose neighbourhoods cannot be taken
disjoint. The only possibility is that $x$ lies on the frontier of $N_1$ in
$M_1$ and $y$ lies on the frontier of $N_2$. Since $N_i$ is GH it is a general
fact that its boundary in $M_i$ is achronal (see~\cite{penrose}).

Thus, up to reverting
time-orientation we may suppose that small timelike curves starting from $x$
must be contained in $N_1$. Denote by $I^{+}_{\epsilon}(x)=\{\exp_x tv|
t\in(0,\epsilon), v\textrm{ future-pointing unit timelike vector at }x\}$ the
set of geodesics of length less than $\epsilon$ starting from $x$. For $\epsilon$
small, $I^{+}_\epsilon(x)\subset N_1$, thus $I^{+}_\epsilon(x)$ isometrically
embedds in $M_2$.
 
Since every neighbourood of $x$ meets every neighbourhood of $y$ in $\hat M$,
and since $\lim_{t\rightarrow 0}\exp_x tv=x$, we have that the image of
$\exp_x tv$ in $M_2$ goes to $y$ as $t\rightarrow 0$. Thus we get that the
identification between $N_1$ with $N_2$ sends $I^+_\epsilon(x)$
to $I^+_\epsilon(y)$. In particular such identification continuously extends at
$x$, sending $x$ to $y$. Now, since $M_1$ and $M_2$ are AdS, we can choose
small neighbourhoods $U,V$ of $x$ and $y$ respectively, such that the isometry
between $I^+_\epsilon(x)$ and $I^+_\epsilon(y)$ extends to an isometry between 
$I^+_\epsilon(x)\cup U$ and $I^+_\epsilon(y)\cup V$. 

This induces an isometry between $N_1\cup U$ and $N_2\cup V$. Up to taking
smaller $U$ and $V$, the spaces $N_1\cup U$ and $N_2\cup V$ are GH, but this
contradicts the maximality of $(N_1, N_2)$.
\end{proof}

In this paper we will deal with convex globally hyperbolic AdS structures with
particles on a fixed topological support $\Sigma\times\mathbb R$ and fixed singular set
equal to $\{p_1,\ldots,p_n\}\times\mathbb R$ (where $p_i$'s are fixed points on
$\Sigma$) with cone angles $\theta_1,\ldots,\theta_n$ respectively. Thanks to
Proposition~\ref{ext:prop} we can restrict to \emph{maximal} convex GH AdS spacetimes.
The corresponding ``Teichm\"uller space'' --- that is the space of GHM AdS structures
on $\Sigma\times\mathbb R$  up to diffeomorphisms isotopic to the identity ---
will be  $\cGH_{\Sigma,n,\theta}$.

\begin{remark}
The condition that the Cauchy surface is locally convex is
perhaps not necessary; it is conceivable that all maximal globally 
hyperbolic AdS manifolds --- containing a Cauchy surface which is perhaps
not locally convex --- actually contain another Cauchy surface which is
locallly convex. We do not elaborate on this point here and simply assume that
this condition is satisfied, since it is necessary in the sequel.
\end{remark}

\subsection*{The left and right representations of a GHMC AdS manifold.}

We have already mentioned in the introduction that $SO_0(2,2)$
is the product of two copies of $PSL(2,\R)$. It follows that the 
representation $\phi_M:\Gamma\rightarrow SO_0(2,2)$ can be identified
with a pair of representations $(\phi_l, \phi_r)$ of $\Gamma$
in $PSL(2, \R)$. 

\begin{lemma} \label{lm:left-right}
$\phi_l$ and $\phi_r$ are the holonomy representations of $\mu_l$ and
$\mu_r$, respectively.  
\end{lemma}

The proof can be found in \cite{minsurf}.

\section{Earthquakes on surfaces with particles}\label{sec:earth}

Let $S$ denote a hyperbolic structure on $\Sigma$ with cone singularities at
$x_1,\ldots, x_n$, with cone angles $\theta=(\theta_1,\ldots,\theta_n)\in
(0,\pi)$.  

If $\lambda$ is a weighted multicurve on $S$  then \emph{the right earthquake} along
$\lambda$ is the fractional negative Dehn twist along each curve with shear factor equal to the
corresponding weight. The corresponding point in $\cT_{\Sigma,n,\theta}$ will
be denoted by $E^r_\lambda(S)$.

The following proposition ensures that the definition of
$E^r_\lambda(S)$ can be extended by continuity to every measured geodesic
lamination. 

\begin{prop}\label{prop:earth:cont}
Let $(\lambda_k)$ be  a sequence of weighted multicurves converging to a
measured 
geodesic lamination $\lambda$. Then the sequence $E^r_{\lambda_k}(S)$ of
hyperbolic surfaces with cone singularities is convergent in $\cT_{\Sigma,
  n,\theta}$.
\end{prop}

Given a measured geodesic lamination $\lambda$, the surface $E^r_\lambda(S)$ is
the limit of $E^r_{\lambda_k}(S)$, where $\lambda_k$ is any sequence of
weighted multicurve converging to $\lambda$.

\begin{cor}\label{cor:earth:cont}
The map
$$ 
\begin{array}{rccc}
E^r_\cdot(S): & \cML_{\Sigma,n} & \rightarrow & \cT_{\Sigma,n,\theta} \\
& \lambda & \mapsto & E^r_\lambda(S)
\end{array}
$$
is continuous.
\end{cor}

The remaining part of this Section will be devoted to prove
Proposition~\ref{prop:earth:cont}.

Denote by $\Sigmat$ the universal covering of
$\Sigma\setminus\{x_1,\ldots, x_n\}$ and by $dev:\Sigmat\rightarrow\mathbb
H^2$ a developing map of $S$.
Given a weighted multicurve $\lambda$, 
let $\lambdat$ be its lifting on $\Sigmat$.
Given an oriented arc $c$ in $\Sigmat$ transverse to $\lambdat$, 
consider the leaves of
$\lambdat$, say $l_1,\ldots, l_k$, cutting $c$.
An orientation is induced on each $l_i$, 
by requiring that at the intersection point with $c$  a positive tangent vector
of $c$ and a positive tangent vector of $l_i$ form a positive basis.

For each $i$, denote by $u(l_i,a_i)$ the element of $PSL(2,\mathbb R)$, representing
a negative translation along $dev(l_i)$ with translation length equal to the
weight $a_i$ of $l_i$. A simple argument shows that the element
\[
     u(l_1,a_1)\circ\cdots\circ u(l_k,a_k)
\]
actually depend only on the endpoints $x,y$ of $c$ and will 
be denoted $\beta_\lambda(x,y)$.

In such a way a map
\[
   \beta_\lambda:\Sigmat\times\Sigmat\rightarrow PSL(2,\mathbb
   R)
\]
is defined. Such a map is a $\pi_1(\Sigma\setminus\{x_1,\ldots,x_n\})$-invariant 
$PSL(2,\mathbb R)$-valued cocycle, that is
\[
\begin{array}{ll}
1. & \beta_\lambda(x,y)\circ\beta_\lambda(y,z)=\beta_\lambda(x,z)\\
2. & \beta_\lambda(\gamma x,\gamma y)=h(\gamma)\beta_\lambda(x,y)h(\gamma)^{-1}
\end{array}
\]
where $h:\pi_1(\Sigma\setminus\{x_1,\ldots,x_n\})\rightarrow PSL(2,\mathbb R)$
is the holonomy representation of $S$.

Moreover, given  a point $x_0\in\Sigmat\setminus\lambdat$ 
the map
$$
\begin{array}{rccc}
dev_\lambda: & \Sigmat & \rightarrow \mathbb
   H^2 \\
& x & \mapsto & \beta_\lambda(x,x_0)dev(x)
\end{array}
$$
is a developing map for $E^r_\lambda(S)$ and 
$$
\begin{array}{rccc}
h_\lambda: & \pi_1(\Sigma\setminus\{x_1,\ldots,x_n\}) & \rightarrow &
PSL(2,\mathbb R) \\
& \gamma & \mapsto & \beta_\lambda(\gamma
   x_0, x_0)\circ h(\gamma) 
\end{array}
$$
is the corresponding holonomy representation.

To prove Proposition \ref{prop:earth:cont} 
it is sufficient to show that  for a sequence
$\lambda_n$ of weighted multicurves converging to a measured geodesic
lamination $\lambda$, the sequence $dev_{\lambda_n}$ converges to a developing
map. This fact is an easy consequence of the following lemma.

\begin{lemma}\label{lem:earth:cont}
If $\lambda_k$ is a  sequence of weighted multicurves converging to a measured geodesic lamination $\lambda$ then
$\beta_{\lambda_k}$ converges to a $\pi_1(\Sigma\setminus\{x_1,\ldots,x_n\})$-invariant $PSL(2,\mathbb
R)$-valued cocycle.
\end{lemma}
\begin{proof}
Since $\beta_{\lambda_k}$ are cocycles it is sufficient to prove that 
$\beta_{\lambda_k}(x,y)$ converges if a geodesic segment $c$ joins $x$ to $y$.

The proof is based on the following lemma. 

\begin{lemma}\label{lem:earth:ort}
Let $\lambda$ be a measured geodesic lamination.
Let $c$ be a geodesic arc in $S$ joining two leaves $l_1,l_2$ of
$\lambda$. Then  either there exists an isometric immersion of a hyperbolic triangle with an ideal vertex
in $S$, sending the compact edge on $c$ and the ideal edges on $l_1$ and $l_2$ respectively or
there exists a geodesic arc $c'$ that satisfies the
following properties:

1. It is homotopic to $c$ through a family of arcs joining $l_1$ to $l_2$ in
$S\setminus\{x_1,\ldots,x_n\}$.

2. It is orthogonal to both $l_1$ and $l_2$.
\end{lemma}

The proof of this lemma will be postponed until the end of this proof.
An easy consequence of Lemma~\ref{lem:earth:ort} is that two
leaves $l$ and $l'$ of $\lambda$ cutting $c$ are sent by $dev$ to
disjoint geodesics.

Take a partition of $c$ in segments $c_1,\ldots,
c_N$ with endpoints $x_i, x_{i+1}$ 
such that the length of each $c_i$ is less than $\varepsilon$. 
For every $k>0$ and  $i=1,\ldots,N$ denote by
$u_{k,i}=\beta_{\lambda_k}(x_i,x_{i+1})$,
$m_{k,i}$ the mass of $c_i$ with respect to $\lambda_k$ and $l_{k,i}$ a leaf
of $\lambda_k$ meeting $c_i$.
If $\varepsilon$ is sufficientely small there exists a constant $C$ such that
\[
   ||u_{k,i}-u(l_{k,i}, m_{k,i})||\leq C\varepsilon m_{k,i}
\]
(the norm we consider is the operatorial norm of $PSL(2,\mathbb R)$
and the inequality is a consequence of the fact that $dev$ sends leaves of
$\lambdat$ cutting $c$ to disjoint geodesics, see Chapter 3 of
\cite{epstein-marden}).

If $c_i$ intersects $\lambdat$ then choose a leaf $l_i$ of
$\lambdat$ and $l_{k,i}$ can be chosen converging to $l_i$ as
$k\rightarrow+\infty$. It follows that up to changing $C$
\[
   ||u_{k,i}-u(l_i,m_{k,i})||\leq C\varepsilon m_{k,i}\,.
\]
Since $\beta_{\lambda_k}(x,y)=u_{k,1}\circ\cdots\circ u_{k,N}$ and $u_{k,i}$'s
 runs in a compact set of $PSL(2,\mathbb R)$, there exists a constant $C'$
 such that
\[
   ||\beta_{\lambda_k}(x,y)-\beta_{\lambda_h}(x,y)||\leq
     C'\sum_{i=1}^N||u_{k,i}-u_{h,i}||\,.
\]
Now if $c_i$ intersects $\lambdat$ then 
\[
   ||u_{k,i}-u_{h,i}||\leq C\varepsilon (m_{k,i}+m_{h,i})+  |u(l_i,m_{k,i})-u(l_i,m_{h,i})|\leq
   C\varepsilon(m_{k_i}+m_{h,i}) + C''|m_{k,i}-m_{h,i}|
\]
for some constant $C''$. Otherwise
\[
    ||u_{k,i}-u_{h,i}||\leq C'''(m_{k,i}+m_{h,i})\,.
\]
for some $C'''>0$.
If $A$ denotes the union of $c_i$'s that does not meet $\lambdat$ then
we get
\[
   ||\beta_{\lambda_k}(x,y)-\beta_{\lambda_h}(x,y)||\leq
   K(\varepsilon (\lambdat_k(c)+\lambdat_h(c)) + \sum_{i=1}^N
   |m_{k,i}-m_{h,i}|+ \lambdat_k(A)+\lambdat_h(A))
\]
for some $K>0$.
Since $\lambdat_k\rightarrow\lambdat$ as $k\rightarrow +\infty$ and
since the mass of $A$ with respect to $\lambdat$ is $0$ by definition, it
follows that $\beta_{\lambda_k}(x,y)$  is a Cauchy sequence.
\end{proof}

\begin{figure}
\begin{center}
\input earth.pstex_t
\end{center}
\end{figure}

\begin{proof}[Proof of Lemma~\ref{lem:earth:ort}]
Consider the family $\mathcal S$ of arcs in $S\setminus\{x_1,\ldots,x_n\}$
homotopic to $c$ through arcs joining $l_1$ to $l_2$ and avoiding the
singularities.

Consider a sequence $c_k\in\mathcal S$ minimizing the length.
Up to passing to a subsequence,
$c_k$ converges either to a point $\hat p$ or to an arc $c_\infty$ in $S$. 

First consider the case that $c_k$ converges to $\hat p$. 
Notice that either $c\cup c_k\cup l_1\cup l_2$ bounds a hyperbolic quadrilateral, 
or they bounds two hyperbolic triangles. In the latter situation the length of $c_k$ 
would be bounded by the length of $c$ multiplied by some constant depending only by 
the angles that $c$ forms with $l_1$ and $l_2$. Thus we can suppose that for all $k$
there exists an isometric immersion of a  hyperbolic quadrilateral 
$Q_k$ in $S$ such that two opposite edges are sent respectively to $c$ and $c_k$ 
and the other two edges, say $u_k,v_k$ are sent to $l_1$ and
$l_2$.  Let $a_k$ and $b_k$ be the lengths of $u_k$ and $v_k$. 
Denote by $p$ and $q$ the end-points of $c$ and by $v,w$ the vectors tangent to $l_1$ 
and $l_2$ respectively at $p$ and $q$. We have that $\exp_p a_kv$ and $\exp_p b_k w$ 
are connected by geodesics shorter and shorter. So, if both $a_k$ and $b_k$ remained 
finite, we whould find an intersection 
point between $l_1$ and $l_2$. 

In the Poicar\'e model we can choose $Q_k$ in such away that the vertex sent to $p$ 
is $0$ and the edges sent to $c$ is contained 
in a fixed geodesic, $\hat c$. Notice that there are two geodesics, say 
$\hat l_1,\hat l_2$ of $\mathbb H^2$ containing all $u_k$ and $v_k$ respectively. 
Thus they cannot intersect  (otherwise $a_k$ and $b_k$ should be bounded) and cannot 
be ultraparallel (otherwise the length of $c_k$ cannot go to $0$). Thus they are 
parallel, and there is a triangle $T$ with an ideal vertex bounded by 
$\hat l_1,\hat l_2$ and $\hat c$. Moreover, $T$ is the union of $Q_k$. Since the 
immersions of $Q_k$'s into $F$ 
coincide on their intersections, it is possible to define an isometric immersion 
of $T$ in $S$ as stated in the Lemma.

Consider now the case $c_k$ goes to an arc $c_\infty$.
If $c_\infty$ avoids
the singularities 
then it is clear that it belongs to $\mathcal S$. Moreover, by standard
variational arguments, it is geodesic and orthogonal to both $l_1$ and $l_2$.

Suppose by contradiction that $c_\infty$ intersects some singuarity. Then it
is not 
difficult to see that it is piecewise geodesic with vertices at singular
points. Since it is homotopic to $c$ there exists an embedded hyperbolic
polygon $P$ whose boundary is contained in $c\cup c'\cup l_1\cup
l_2$. Moreover, at least 
one vertex of this polygon can be supposed to lie on a singular point. 
Since the cone angles are less than $\pi$, $P$ is convex. It follows that the
length of $c'$ can be shortened.
\end{proof}

 
\section{Convex subsets in AdS manifolds with particles}

\subsection*{Globally hyperbolic convex subsets.}

Let $M$ be an AdS manifold with particles. We are interested in convex subsets
in $M$, and in particular in the ``convex core'' which will be defined
below. This convex core always contains a closed space-like surface, and 
this is also the case of many (perhaps all) convex subsets. For technical
reasons, we are lead to include this property in the definition of convex
subsets, so that we consider here ``globally hyperbolic'' convex subsets.

\begin{defi}
Let $\Omega$ be a non-empty, connected subset of $M$. It is {\it GH convex} if:
\begin{itemize}
\item $\Omega$ contains a closed, space-like surface $S$,
\item $\Omega$ has a space-like, locally convex boundary. 
\end{itemize}
\end{defi}

It follows from the considerations made below that this definition is
equivalent to a ``global'' definition: any geodesic segment in $M$, with
endpoints in $\Omega$, is actually contained in $\Omega$. However the
definition given here is more convenient here.
We will often write ``convex'' instead of
``GH convex'', hoping that it does not confuse the reader.

\subsection*{Convex GHM AdS manifolds contain a compact GH convex subset.}

The reason why we consider only {\it convex} GHM manifolds is that they always
contain a compact, GH convex subset. This will be used below to show that they
actually contain a smallest such convex subset, their ``convex core''.

\begin{lemma} \label{lm:52}
Suppose that $M$ contains a closed, space-like, locally convex surface
$S_0$. Then $M$ contains another such surface, arbitrarily close to $S_0$,
which in addition is smooth.  
\end{lemma}

For the proof of this lemma, we will use the notion of the distance from
$S_0$. It is defined at $p$ as the \emph{maximum} of the Lorentzian length
of causal segments connecting $p$ to $S_0$.

\begin{lemma}\label{dist:lem}
Let $S$ be a convex spacelike surface in $M$, and let $d$ be the distance from
$S$. Then
\begin{enumerate}
\item $d$ is continuous;
\item for every point $p$ there is a geodesic segment joining $p$ to $S$ that
  avoids the singularities and realizes the distance;
\item if $d(p)<\pi/2$ and $p$ is in the convex region bounded by $S$, then
  this segment is unique;
\item the function $d$ is $\mathrm C^{1,1}$ on the set of points in the
  bounded region by $S$ satisfying $d(p)<\pi/2$;
\item if $S$ is in the class $\mathrm C^k$, then $\mathrm C^{1,1}$ can be
  replaced by $\mathrm C^{k}$ in the previous point.
\end{enumerate}
\end{lemma}

\begin{proof}
The continuity of $d$ is an easy consequence of the compactness of $S$. 

For point (2), it is well known that 
the space of causal curves joining $p$ to $S$ is compact and
the length function is upper-semicontinuous (for a proof see Sections 6 and 7
of~\cite{penrose}). Thus a causal path realizing the distance exists. A priori
this path is piece-wise geodesic with vertices on the singular locus. 
Suppose that a vertex, say $r$, on a cone singularity of angle $\theta_i$, 
occurs. 
Then take points on the curve, $r_-, r_+$ respectively in the past 
and the future of $r$.
If $r_-,r_+$ are chosen in some neighbourhood of  $r$ isometric to 
a convex neighbourhood 
of $H^3_{1,\theta_i}$, then the segment joining $r_-$ and $r_+$ 
has length bigger 
than the sum of the lengths joining $r_-$ to $r$ and $r$ to $r_+$. 
Thus the curve can be lengthened.

For point (3) we consider the unit normal bundle of the surface $S$,
i.e. the set of unit vectors $n\in T_xS$, for $x\in S$, such that 
the plane orthogonal to $n$ at $x$ is a support plane of $S$ and $n$
is towards the convex side of $S$. This unit normal bundle, $N^1S$,
is a submanifold of the tangent bundle of $M$. The normal exponential
map of $S$ is the map $\exp:N^1S\times \R \rightarrow M$ sending a unit
vector $n\in N^1S$ and a $t\in \R_{>0}$ to the endpoint of the time-like 
geodesic segment of length $t$ starting from the basepoint of $n$ with
initial velocity equal to $n$. We denote this endpoint by $\exp_t(n)$.

Let $t\in (0,\pi/2)$, and let $S_t:=\exp_t(N^1S)$. Clearly
any point on the convex side of $S$ which can be connected to $S$ by
a maximizing (time-like) geodesic segment $g$ of length $t$ has to be in 
$S_t$, since the plane orthogonal to $g$ at its intersection with $S$
is a support plane of $S$ (my maximality of $g$). 

We claim that $\exp_t$ is a homeomorphism from $N^1S$ to $S_t$. The 
fact that $\exp_t$ is a local homeomorphism follows from the properties
of Jacobi fields along time-like geodesics in $AdS^3$. It is not 
difficult to check that such Jacobi fields are of the form
$$ J(s) = \cos(s) v_0 + \sin(s) v_1~, $$
where $v_0$ and $v_1$ are parallel vector fields 
along the time-like geodesic.
Given $n\in N^1S$, a first-order displacement of $n$ on $N^1S$
induces a Jacobi field $J(s)$ along the geodesic segment with initial
velocity equal to $n$ which is of the form above, with 
$\langle v_0, v_1\rangle \geq 0$ -- because $S$ is convex -- but
either $v_0$ or $v_1$ non-zero. If $v_0\neq 0$, then 
$$ \langle J(s), v_0\rangle = \cos(s) \| v_0\|^2 + \sin(s) 
\langle v_0,v_1\rangle >0 $$
for all $s\in (0,\pi/2)$, while if $v_1\neq 0$ then 
$\langle J(s), v_1\rangle>0$ for all in $s\in (0,\pi/2)$. 
Therefore 
$J(s)$ does not vanish for $s\in (0,\pi/2)$, which means that
$\exp_t$ is a local homeomorphism from $N^1S$ to $S_t$ for 
$t\in (0,\pi/2)$.

Suppose that $\exp_t$ is not a global homeomorphism for some $t\in 
(0,\pi/2)$. Let $t_0$ be the infimum of all $t\in (0,\pi/2)$ for
which $\exp_t$ is not a global homeomorphism, then there are points
$n_1, n_2\in N^1S$ such that $\exp_{t_0}(n_1)=\exp_{t_0}(n_2)$. It follows
from the definition of $t_0$ that the geodesic segments 
$\exp_{[0,t_0]}(n_1)$ and $\exp_{[0,t_0]}(n_2)$ are parallel at their
common endpoint, this is clearly impossible unless $n_1=n_2$. So $\exp_t$ is a 
global homeomorphism from $N^1S$ to $S_t$ for all $t\in (0,\pi/2)$,
and this proves point (3).

For point (4), because of point (3)
one can adapt to the AdS context a nice argument developed for the
corresponding hyperbolic situation by Bowditch (see \cite{epstein-marden} and
\cite{mess,be-bo} for the Lorentzian case).
\end{proof}

\begin{proof}[Proof of Lemma \ref{lm:52}]
For each $\pi/2>\epsilon>0$, let $S_\epsilon$ be the set of points at distance 
$\epsilon$ from $S_0$ on its concave side. 

Then $S_\epsilon$ is not only convex, but it is uniformly convex: given a
point $p\in S_\epsilon$ let $q$ be a point on $S$ such that the distance
between $p,q$ is $\epsilon$ (in general such a point is not unique). The germ
of plane $P$ through $q$ and orthogonal to the geodesic segment between $p$
and $q$, turns out to be a support plane for $S$. Moreover the set $P_\epsilon$ of
points at distance $\epsilon$ from $P$ is tangent to $S_\epsilon$ at $p$ and
in fact $S_\epsilon$ is contained in the convex region bounded by $P_\epsilon$
in a neighbourhood of $p$.  The uniform convexity of $S_\epsilon$ follows
then since the set of points at distance $\epsilon$ from a geodesic plane in
$AdS_3$ is uniformly strictly convex (actually even umbilic). However
$S_\epsilon$ is not smooth.

Now choose $\epsilon'\in (0,\epsilon)$ small enough, and let 
$S'_{\epsilon,\epsilon'}$ be the surface at constant distance
$\epsilon'$ from $S_\epsilon$ on the {\it convex} side. Note that
$S'_{\epsilon,\epsilon'}$ is {\it not} the surface at constant 
distance $\epsilon-\epsilon'$ from $S$. 

Notice that not only $S'_{\epsilon, \epsilon'}$
converges to $S_\epsilon$ as $\epsilon'\rightarrow 0$, but also the unit
normal bundle of $S'_{\epsilon, \epsilon'}$ converges to the unit normal
bundle of $S_\epsilon$ in $TM$. Thus the unifom convexity of
$S_\epsilon$ easily implies that $S'_{\epsilon,\epsilon'}$ is convex (and also
uniformly convex) for $\epsilon'$ close to $0$.

The last step is to smooth $S'_{\epsilon,\epsilon'}$ to obtain a surface $S$
as needed. Standard arguments (based on convolution and partitions of
unity) can be used here to obtain a smooth embedding which is $C^1$-close
to $S'_{\epsilon,\epsilon'}$; since the question is of a local nature, 
the corresponding (classical) results in the Euclidean 3-space can 
actually be used here through the projective model of AdS.
\end{proof}

\begin{lemma} \label{lm:existence}
Let $M$ be a convex GHM AdS manifold with particles. Then $M$ contains a
compact, GH convex subset.
\end{lemma}

\begin{proof}
By definition of a convex GHM manifold, $M$ contains a closed, space-like,
locally convex surface $S_0$ which is orthogonal to the particles. 
By the previous remark $M$ also contains a closed, space-like, locally convex
surface $S$, which in addition is smooth.

Let $I$ and $B$ be the induced metric and shape operator of $S$,
respectively. Consider the manifold $S\times [0,\pi/2)$ with the
metric: 
$$ g_S := 
-dt^2 + I((\cos(t)E+\sin(t)B)\cdot, (\cos(t)E+\sin(t)B)\cdot)~, $$
where $t$ is the coordinate in $[0,\pi/2)$ and $E$ is the identity
morphism on $TS$. A simple computation shows that the metric 
$g_S$ is locally modeled on $H^3_1$ --- except at the points
which project to singular points of $S$, where there are of course
cone singularities. Also by construction (and because $M$ is
maximal), $(S \times [0,\pi/2), g_S)$ embeds isometrically into 
$M$, with the surface $S\times \{ t\}$ sent to the set of points
at distance $t$ from $S$ on the convex side. Finally note that
$(S \times [0,\pi/2), g_S)$ has locally convex boundary,
and that its shape operator at $t=\pi/2$ is simply $B^{-1}$ (unless
$B$ is positive semi-definite but not positive definite, in which
case the boundary has ``pleating lines''). 
So the closure of the image of $S \times [0,\pi/2)$ in $M$ is a GH 
convex subset which is compact, as required.
\end{proof}

\subsection*{The distance to a convex subset is bounded.}

Another key property of GHM AdS manifolds is that the distance from a GH
convex subset to the boundary is always less than $\pi/2$. 

\begin{lemma}\label{lm:distance}
Let $M$ be a convex GHM AdS manifold, and let $K$ be a GH convex, compact
subset of $M$. For each $x\in M\setminus K$, the maximal time-like geodesic
segment(s) joining $x$ to $K$ has (have) length less than $\pi/2$. 
\end{lemma}

The proof is based on a simple proposition concerning the closest point
projection on a convex subset in an AdS manifold. We call $\exp_K$ the normal
exponential map, which is defined on the unit normal bundle of $\dr K$, 
and $\exp_K^r$ the restriction to the set of vectors of norm equal to $r$.

\begin{prop} \label{pr:3pts}
Let $M$ be a GHM AdS manifold, and let $K$ be a GH convex subset of $M$.
Let $x\in M\setminus K$ be a regular point of $M$ at distance less than 
$\pi/2$ from $K$. 
\begin{enumerate}
\item There exists a (time-like) 
geodesic segment $\gamma$ which has maximal length going from $x$ to 
$\dr K$, which connects $x$ to a point $y\in \dr K$. Let $r$ be its length.
\item $\dr K$ is $C^{1,1}$ smooth at $y$, with principal curvatures bounded
  from below, and its shape operator $B$ is such
  that  $\cos(r) E -\sin(r) B$ is non-negative.
\item Let $P$ be the support plane of $K$ at $y$, let $\Pi':T_yM\rightarrow
T_xM$ be the parallel transport along $\gamma$, and let $\Pi$ be its
restriction from $P$ to its image. Then the map $\Pi^{-1}\circ
d\exp_K^r:P\rightarrow P$ is equal to:
$$ \Pi^{-1}\circ d\exp_K^r = \cos(r) E -\sin(r) B~. $$
\end{enumerate}
\end{prop}

\begin{proof}
We suppose without loss of generality that $x$ is in the future of $K$.
Since $M$ is GH and $K$ contains a Cauchy surface, any past-oriented
time-like curve starting from $x$ intersects $K$. In particular this
holds for all time-like geodesics, so a simple compactness 
argument based on the time cone of $x$ shows that there exists
a time-like geodesic segment of maximal length connecting $x$ to $K$.
This segment is not necessarily unique.

For point (2) note first that, since $\gamma$ has maximal length, the plane
$P$ orthogonal to $\gamma$ at $y$ is a support plane of $K$. Fix a 
small neighborhood $U$ of $y$, and let $H$ be the set of points in $U$ which
can be joined to $x$ by a time-like geodesic segment of length $r$.
Since $r\in (0,\pi/2)$, $H$ is a space-like umbilic surface of principal
curvatures equal to $\cotan(r)$. The definition of $y$ shows that 
$H$ is tangent to $P$ at $y$. But $y$ maximizes the distance from
$x$ to $\dr K$, it follows that, in $U$, $\dr K$ is in the future of 
$H$ (because any geodesic segment going from $x$ to $H$ has to intersect
$\dr K$ before $H$). So $\dr K$ is ``pinched'' at $y$ between $P$ and $H$,
which implies that it is $C^{1,1}$ with principal curvatures bounded 
between $0$ --- the principal curvatures of $P$ --- and $\cotan(r)$ --- the
principal curvatures of $H$. This proves point (2).

For point (3) consider a geodesic segment $\gamma:[0,r]\rightarrow AdS_3$, and
let $u$ be a unit vector field along $\gamma([0,r])$ which is parallel and
orthogonal to $\gamma([0,r])$. For all $v,w\in \R$ there is a unique Jacobi
field $Y$ along $\gamma([0,r])$ with $Y(0)=vu$ and $Y'(0)=wu$, it is equal to
$Y(s)=(v\cos(s)+w\sin(s))u$ at $\gamma(r), 0\leq s\leq r$. We can apply this
computation with $u$ equal to a principal vector at $x$, $v=1$, and $w$ equal
to the principal curvature corresponding to $x$, this yields point (3).
\end{proof}

\begin{proof}[Proof of Lemma \ref{lm:distance}]
We suppose again, still without loss of generality, that $x$ is in the future
of $K$. 

The first point is that there exists a space-like curve $c:[0,L]\rightarrow M$
which begins on $\dr K$ and ends at $x$. This follows again from the global
hyperbolicity of $M$; since every past-directed time-like and light-like
curve starting from $x$ intersects $K$, it is also true for some space-like
curves starting from $x$. 

Consider the distance to $K$, it is a function, which we still call $r$,
defined in $M\setminus K$. We suppose (by contradiction) that $x=c(L)$
is at distance at least $\pi/2$ from $K$, and let $t_0$ be the minimum of all
$t\in [0,L]$ for which $c(t)$ is at distance at least $\pi/2$ from $K$.

The function $r$ is continuous but not $C^1$, since there
are points which are joined to $K$ by two maximizing (time-like) geodesic
segments. Notice that the function $r$ could be regarded as the cosmological 
time of the spacetime $I^+(K)$ -- that is $r(x)$ is the $\sup$ 
of the Lorentzian lengths of the timelike curves contained in $I^+(K)$ with future end-point at $x$.
Thus , from the general result in~\cite{AGH}, 
$r$ is a 
semi-convex function, i.e. in local charts it
is the sum of a convex function and a smooth function. In particular $r$ 
is twice differentiable almost everywhere. 
It follows that it is possible to choose a generic  space-like curve $c$, which
intersects the points where $r$ is not differentiable on a set $E$ of measure $0$,
so that $r$ is continuous and, for $s\not\in E$, $r$ is differentiable at $s$  
and there is a unique maximizing segment between $c(s)$ and $K$.
At such points, we call $\pi(s)$ the endpoint on $\dr K$ of the 
maximizing curve between $x$ and $K$. 
Let us stress that, since $r$ is semi-convex, then $(r\circ c)'\in L^\infty([0,L])$, 
that implies that $r\circ c$ is an absolutely continuous function (it coincides 
with the integral of its derivative).

So $\pi:[0,t_0]\rightarrow \dr K$ is well-defined and continuous 
except on $E$, and it is Lipschitz on the complement of $E$; its derivative
$\pi'(s)$ is defined almost everywhere. We now suppose that 
$c$ is parametrized in such a way that $\pi$ is parametrized
at speed $1$ on the complement of $E$.

It follows from point (3) of the previous proposition that for 
$s\not\in E$ the norm of the
image of $\pi'(s)$ by the map $\exp_K^{r(s)}$ is equal to 
$\| (\cos(r)E-\sin(r)B)\pi'(s)\|$, so that (using point (2)
of the proposition) it is bounded by $\cos(r)$. But this vector 
is the ``horizontal'' component of $c'(s)$ (its projection on the
kernel of $dr$ in $T_{c(s)}M$). Since $c$
is space-like, it follows that $|(r\circ c)'(s)|\leq \cos(r)$.
So the function $r\circ c$ is well-defined on $[0,t_0]$,
absolutely continuous, and it is a solution outside $E$
of the differential inequality: $y'(s)\leq \cos(y(s))$, with
$y(0)=0$ by definition of $c$. It follows quite directly  that
$y(t_0)<\pi/2$, which contradicts the definition of $t_0$. This shows
that $x$ is at distance less than $\pi/2$ from $K$, as announced.
\end{proof}

\subsection*{The existence of the convex core.}

It is now possible, using in particular Lemma \ref{lm:distance}, to show that
the intersection between two GH convex subsets of $M$ is itself non-empty and
GH convex. This is a key point in proving the existence of a ``convex core''
in $M$. 

We will use the following simple remark.

\begin{lemma}\label{bound:lem}
Let $\Omega$ be a GH convex subset of an AdS spacetime $M$.
Each boundary component of $\Omega$ is either convex in the past or convex in the future.
(A locally convex spacelike surface $S$ is said to be convex in the future
(resp. past) if future-oriented timelike unit vector normal to $S$ points towards the convex
(resp. concave) side bounded by $S$.)

Suppose $\Omega$ to be compact. Then there are $2$ boundary components:
the \emph{future boundary}, $\dr_+\Omega$ that is convex in the past, and
the~\emph{past boundary}, $\dr_-\Omega$ that is  convex in the future.

Both $\dr_-\Omega$ and $\dr_+\Omega$ are Cauchy surfaces and
\begin{equation}\label{lll}
   \Omega= I^+(\dr_-(\Omega))\cap I^-(\dr_+(\Omega))
\end{equation}
\end{lemma}

\begin{proof}
For the first part, it is sufficient to notice that if the timelike vectors
tangent to $x\in S$ points towards the convex side bounded by $S$, then the
same holds for points $y$ in a neighbourhood of $x$ in $S$. 

For the second part, the claim is that any inextensible causal curve intersect
$\Omega$ in a compact interval whose future end-point lies on $\dr_+\Omega$
and the past end-point lies on $\dr_-\Omega$.
Since $\Omega$ is compact in $M$ then there is a Cauchy surface, say $S_+$, in the future
of $\Omega$ and a Cauchy surface, say $S_-$ in the past of $\Omega$. The curve
$c$ intersects $S_-$ in a point $x_-$ and $S_+$ in a point $x_+$. It is clear
that $c\cap\Omega$ is contained in the segment bounded by $x_-$ and $x_+$ on $c$.
Now suppose $I$ to be a connected component of $c\cap\Omega$. By definition we
should have that the past endpoint of $I$ lies on $\dr_-\Omega$ whereas the
future endpoint lies on $\dr_+\Omega$. If $c\cap\Omega$ contained two connected
components, then one could construct a  causal curve in $M$ with past
end-point, say $x_P$, on $\dr_+\Omega$ and future end-point on
$\dr_-\Omega$. Gluing $c$ with a timelike geodesic with future end-point at
$x_P$ and a timelike geodesic with past end-point at $x_F$, produces a causal
curve meeting $S$ twice.

Consider now the flow $\phi$ of some future-oriented unit
timelike vector-field $X$ on $M$. It
follows from the claim that the set, $\mathcal S$, of points $(x,t)\in
S\times\mathbb R$ such that
$\varphi_t(x)\in\Omega$ is a compact regular
neighbourhood of $S$, that is $\mathcal S\cong S\times[-1,1]$ in such a way
the image through $\phi$ of $S\times\{1\}$ is $\dr_+ \Omega$ and the image
of $S\times\{-1\}$ is $\dr_-\Omega$.

The identity~(\ref{lll}) is a simple consequence of the claim.
\end{proof}

\begin{prop}
Let $\Omega, \Omega'$ be two GH convex subsets of $M$. Then $\Omega\cap
\Omega'$ is non-empty and GH convex.
\end{prop}

\begin{proof}
We claim that $\dr_+\Omega$ is contained in the future of $\dr_-\Omega'$
Suppose a point $x\in\dr_+\Omega$ is contained in the past of
$\dr_-\Omega'$. Let $M'$ be the AdS structure on $\dr_+\Omega\times[0,\pi/2)$ 
with metric given by
\[
-dt^2 + I((\cos(t) E+\sin(t) B)\cdot, (\cos(t) E + \sin(t) B)\cdot)
\]
where $I$ is the first fundamental form on $\dr_+\Omega$, $B$ is the shape
operator, and $t\in [0,\pi/2)$. 
The past of $\dr_+\Omega$ is isometric to a regular neighbourhood of
$S\times\{0\}$ in $M'$. If we glue the future of $\dr_+\Omega$ to $M'$, we
get a smooth spacetime $M''$ that contains $M$. The surface $\dr_+\Omega$
turns to be a convex Cauchy surface so $\dr_+\Omega$ is convex GH.

Now in $M''$ there is a timelike geodesic with length equal to $\pi/2$ and
future end-point at $x$. Since $x$ is in the past of $\dr_-\Omega'$, there is a
timelike curve in the complement of $\Omega'$ with length bigger than
$\pi/2$. This contradicts Lemma~\ref{lm:distance}.

The fact that $\Omega\cap\Omega'$ is a non-empty convex subset follows
directly from  the claim and from~(\ref{lll}).
Moreover, the claim implies that every inextensible causal curve meets
$\Omega\cap\Omega'$ in a non-empty interval.
So the same arguments  as in Lemma~\ref{bound:lem} show that $\Omega\cap\Omega'\cong
S\times [0,1]$ and the boundary components of $\Omega\cap\Omega'$ are Cauchy
surfaces.
\end{proof}

It follows from the previous lemma and from Lemma \ref{lm:existence} 
that $M$ contains a minimal GH convex subset.

\begin{lemma}
Let $M$ be a GHMC AdS manifold, then $M$ contains a non-empty GH convex subset
$C(M)$ which is minimal: any non-empty GH convex subset $\Omega$ in $M$
contains $C(M)$.
\end{lemma}

\begin{proof}
We already know that $M$ contains a 
GH convex subset (Lemma \ref{lm:existence})
and that the intersection of two GH convex subsets is GH convex.
We can therefore consider the intersection $K$ of all GH convex 
subsets in $M$, it is clear that it intersects all time-like curves
in $M$ and that it has locally convex boundary. The only point 
that remains to prove is that it has space-like boundary, since 
a limit of space-like surfaces could a priori be light-like.
We do the proof here for the future boundary of $K$, denoted by
$\dr_+K$, the same argument applies with obvious changes to the
past boundary $\dr_-K$.

Let $S_\pm$ be respectivey the future and the past boundary of
the GH convex compact subset appearing in the proof of Lemma~\ref{lm:existence}.
By construction, the distance between $S_+$ and $S_-$ is $\pi/2$.
Since $K$ is contained in any GH convex sets,  all
points of $\dr_+K$ are at distance less than $\pi/2$ from $S_-$.

Let $\Omega$ be the set of points at distance less
than $\pi/2$ in the future of $S_-$. 
We consider the function $u$ defined, on $\Omega$, 
as the sine of the distance to $S_-$. 
By Lemma~\ref{dist:lem} $u$  is a smooth function on $\Omega$. Moreover it satisfies the
equation:
$$ \hess(u) \leq ug~, $$
where $g$ is the AdS metric on $M$. To check this equation note that it is
satisfied (and is actually an equality) if $S_-$ is totally geodesic. If $x\in
\Omega$ and if $v\in T_x\Omega$ is the direction of the maximal geodesic
segment from $x$ to $S_-$, the Hessian of $u$ behaves on $\R v$ as the Hessian
of the distance to a geodesic plane, and on the plane orthogonal to $v$ it is
smaller since $S_-$ is convex. 

Now suppose that there exists a point $x\in \dr_+K$ where $\dr_+K$ is
light-like. Then, since $\dr_+K$ is convex, there exists a light-like
past-oriented geodesic ray $\gamma$ contained in a support plane of $\dr_+K$
at $x$. Let $(\gamma_k)$ be a sequence of space-like geodesic rays converging
to $\gamma$, parametrized at unit speed. Then 
$$ \lim_{k\rightarrow
\infty}(u\circ \gamma_k)'(0)=-\infty~, $$
and $(u\circ \gamma_k)''\leq u\circ \gamma_k$ by the estimate on the Hessian
of $u$. Therefore, for $k$ large enough, $\gamma_k$ intersects $S_-$
at time $t_k$, with $\lim_{k\rightarrow \infty} t_k=0$. It follows that
$\gamma$ intersects $S_-$, this contradicts the convexity of $\dr_+K$.
\end{proof}

As for non-singular AdS manifolds, we call $C(M)$ the {\it convex core} of
$M$.

\section{Pleated surfaces  in AdS manifolds with particles}

\subsection*{The geometry of the convex core}

The boundary of the convex core of a GHM AdS manifold with particles shares
all the important properties of the boundary of a non-singular GHM AdS
manifold (as studied in \cite{mess}), which are also the same as for 
quasi-Fuchsian hyperbolic manifolds (including those with particles as in
\cite{minsurf,qfmp}).  

The first property is that boundary components of the convex core are pleated
surfaces according to the following definition.

\begin{defi}
A convex {\bf pleated surface} in $M$ is a closed, convex, space-like surface,
$S$ orthogonal to the singular locus of $M$, which is ruled: for any point
$x\in S$, other than in the singular set of $M$, $x$ is contained in 
the interior either of a
geodesic segment of $M$ contained in $S$, or of a geodesic disk contained in
$S$. 
\end{defi}

\begin{lemma}\label{lm:bd}
Let $M$ be a convex GHM AdS manifold with particles.
Each boundary component of its convex core $C(M)$ is a pleated surface.
If $p\in \dr C(M)$ does not lies on a singular line, then $C(M)$ has a
unique support plane at $p$, say $H$, and $H\cap\dr C(M)$ contains a 
neighbourhood of $p$ in $H$. 
\end{lemma}
 
\begin{proof}
Suppose by contradiction that some vertex occurs, that is, there exists a
support plane $P$ intersecting $\dr C(M)$ at exactly one point $p$. Without
loss of generality we may suppose $p\in\dr_+ C(M)$. 

Let $v$ be the unit vector orthogonal to $P$ at $p$ and pointing in its past
and consider the plane 
$Q_\epsilon$ orthogonal to the geodesic $\gamma(t)=\exp_p tv$ at
$\gamma(\epsilon)$.
For $\epsilon$ sufficiently small, 
then $Q_\epsilon\cap\dr_+ C(M)$ is topologically a circle and
$Q\cap\dr_- C(M)=\varnothing$.  Let $\Delta'$ the surface obtained by
replacing  in $\dr_+ C(M)$ the set  $I^+(Q_\epsilon)\cap\dr_+ C(M)$ with $Q_\epsilon
\cap C(M)$. Then $\Delta$ is a locally convex in the past surface. 
The convex domain
$\Omega= I^{-}(\Delta)\cap I^+(\dr_- C(M))$ is then smaller than $C(M)$.

To show that $\dr C(M)$ is orthogonal to the singular locus of $M$, 
let $x\in \dr C(M)$ be a point contained in a singular curve of $M$. 
We consider the link $L_x(M)$ of
$M$ at $x$ (the set of geodesic rays in $M$ starting from $x$, with the
natural angular metric). It is a (real) projective surface with two singular
points. It is also endowed naturally with a ``distance'', coming from the
angles between the geodesic rays starting from $x$, which is locally modeled
on the de Sitter plane for space-like rays, and on the hyperbolic plane for
the time-like rays (except the two rays which follow the singular line of
$M$ containing $x$). $L_x(M)$ also contains a closed curve $\gamma_0$,
which is the union of rays orthogonal to the singular line containing $x$,
it is geodesic for the metric just described (and a line for the real
projective structure on $L_x(M)$). We call $h$ the oriented 
distance to $\gamma_0$
(i.e., the length of the maximal geodesic connecting a point to $\gamma_0$),
with a plus sign if this segment is past-oriented, and a minus sign if it is
future-oriented). 
It is well-defined on the ``de Sitter'' part of $L_x(M)$ (corresponding to the
space-like rays). 

In $L_x(M)$ we consider the link $L_x(C(M))$ 
of $C(M)$ at $x$, namely, the set of geodesic
rays starting from $x$ for which a neighborhood of $x$ is contained in
$C(M)$. Since $\dr C(M)$ is space-like, $\dr L_x(C(M))$ is a space-like curve
contained in the ``de Sitter'' part of $L_x(M)$. 
It follows from the convexity of $C(M)$ that $\dr L_x(C(M))$ is a locally
convex curve. 

\begin{figure}[ht]
\centerline{\psfig{figure=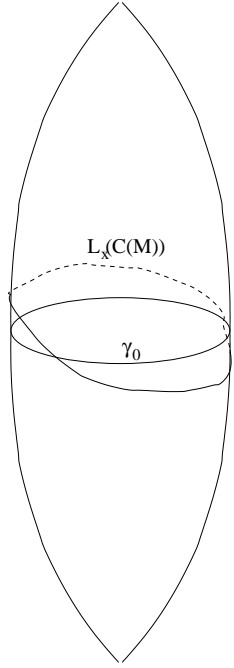,width=3cm}}
\caption{The link of $C(M)$ at a singular point.} 
 \label{fig:link}
\end{figure}

Let $y$ be point of $\dr L_x(C(M))$ where $h$ attains its minimum. 
Then there is a geodesic segment $\gamma_y$ in $L_x(M)$ containing $y$, which
is a support line of $\dr L_x(C(M))$ at $y$, and such that the restriction of
$h$ to $y$ is extremal at $y$. Considering the geometry of $L_x(M)$ shows that
there is a maximal extension of $\gamma_y$ as an embedded curve in $L_x(M)$
which is symmetric with respect to $y$ (it has the same length on both sides
of $y$). 

The fact that the cone angle of $M$ at the line containing $x$ is less than
$\pi$ then shows that the restriction of $h$ to $\gamma_y$ is everywhere
positive, and that it has a maximum at $y$. However the convexity of $\dr
L_x(C(M))$ shows that it is ``below'' $\gamma_y$: any geodesic orthogonal to
$\gamma_0$ intersects $\gamma_y$ and $\dr L_x(C(M))$ once, and the value of
$h$ at the intersection with $\gamma_y$ is bigger than the value of $h$ at the
intersections with $\dr L_x(C(M))$. 

This already shows that $h(y)$ can not be positive, otherwise the restriction
of $h$ to $\dr L_x(C(M))$ would have a strict maximum at $y$, contradicting
the definition of $y$. The same argument shows that if $h(y)=0$ then $h=0$
everywhere on $\dr L_x(C(M))$, i.e., $\dr L_x(C(M))=\gamma_0$.

But $h(y)$ can not be negative, otherwise $h$ would be negative everywhere on
$\gamma_y$, and therefore on $\dr L_x(C(M))$. This would mean that any plane
intersecting the singular line containing $x$ orthogonally a little ``below''
$x$ would cut a small cone off $C(M)$, leaving a piece of $C(M)$ which would
remain GH convex, and would thus contradict the definition of $C(M)$ as a minimal
GH convex set. 

So the only possibility is that $\dr L_x(C(M))=\gamma_0$, so that $\dr C(M)$
is orthogonal at $x$ to the singular locus of $M$. 
\end{proof}

Let the bending locus of $\dr C(M)$, say $L$, be the complement of points that admit
some support plane $P$ such that $P\cap\dr C(M)$ is a neighbourhood of $p$ in
$\dr C(M)$. $L_\pm$ denote the intersection of $L$ with $\dr_\pm C(M)$.

\begin{lemma}\label{lm:fol}
If $L_+=\varnothing$ (resp. $L_-$)
then $C(M)=\dr C_-(M)=\dr C_+(M)$ is a totally  geodesic surface
orthogonal to the singular locus.

If $L\neq\varnothing$ then $L$ is foliated by complete spacelike geodesics of
$M$. 
\end{lemma}

\begin{proof}
If $L_+$ is empty then $\dr_+ C(M)$ is totally geodesic.
In particular it is a convex subset, so it is contained in the convex core.

The second part is more delicate.
Suppose  $p\in L$. There exists a unique $v\in T_pM$ such that the segment
$\gamma$ defined by 
$\gamma(t)=\exp_p{tv}$ is contained in $L$ for $t\in(-\epsilon,\epsilon)$. 
Now consider the set
\[
 A=\{t\in\mathbb R|\exp_p tv\in L\}\,.
\]
Since $L$ is closed in $\dr C(M)$ the set $A$ is closed.
We will show that is open too.
In fact it is sufficient to show that $\exp_p{tv}\in L$ for $t$ small, because
the fact that $A$ contains a neighbordhood of any of its points, say $t_0$, then
follows by substracting $t_0$ to $t$ for $t$ close to $t_0$.

Suppose by contradiction that there is support plane $Q$ that intersects
$\dr_+ C(M)$ in a neighbourhood of $\exp_p tv$. Since $\gamma(s)$ is
contained in $\dr_+ C(M)$ for $s\in (0,t)$ we have that $\gamma'(t)$ is contained in $Q$. It
follows that $\gamma$ is contained in $Q$. Thus $Q\cap \dr_+C(M)$ contains the
convex hull of a ball of center $\gamma(t)$ in $Q$ and $\gamma(-t)$, that is, a
neighbourhood of $p$. This contradicts the assumption that $p\in L$.
\end{proof}

Notice that $M$ contains a geodesic Cauchy surface iff $L=\varnothing$ 
and this is he only case where the convex core has empty interior.
The leaves of the foliation of $L$ pointed out in Lemma \ref{lm:fol} 
will be called the bending lines of the convex core.

\begin{lemma}\label{lm:bd2}
The induced metric on $\dr C(M)$ is hyperbolic,
with cone singularities at the intersections with the particles,
of angle equal to the angle of $M$ at the corresponding singular
lines. The bending locus $L$ is at finite distance
from the cone points. 
\end{lemma}
 
\begin{proof}
A direct consequence of the orthogonality result in Lemma \ref{lm:bd}
is that the induced metric on $\dr
C(M)$ has cone singularities at the intersection of $\dr C(M)$ with the
singular curves of $M$, and that the cone angles at that points are the same
as the cone angles at the corresponding singular lines of $M$.  

This also shows that no leave of the bending lamination of $\dr C(M)$ can go
through $x$, since otherwise $\dr L_x(C(M))$ could not be geodesic. So the
leaves of the bending lamination are embedded geodesics (for the induced
metric on $\dr C(M)$, and the fact that the cone angle at $x$ is strictly less
than $\pi$ then implies that those geodesics can not enter a ball centered at
$x$ (of radius depending on the cone angle at $x$).
\end{proof}

\begin{prop}\label{adstohyp:prop}
$\dr C(M)$ carries an intrinsic $\mathrm C^{0,1}$-hyperbolic structure. 
$L$ is the support of a measured geodesic
lamination $\lambda$ on $\dr C(M)$, called the bending lamination. 
The hyperbolic structure on $\dr_+ C(M)$ and the measured lamination 
$\lambda_+=\lambda|_{\dr_+ C(M)}$ determine $M$.
\end{prop}

This statement is the analog with cone singularities of a well-known fact
for non-singular hyperbolic metrics. We include a proof for the reader's
convenience because the proof in the non-singular case relies heavily on
the use of the developing map and therefore does not extend to hyperbolic
surfaces with cone singularities.

\begin{proof}
First notice that we can choose coordinates around a point of $p\in\dr C(M)$, 
such that $\dr C(M)$ looks like
the boundary of a convex set of $\mathbb R^3$. This show that points of $\dr C(M)$ 
are locally connected by Lipschitz paths.
Thus each component of $\dr C(M)$ is connected by Lipschitz
paths. 
So we can consider the path distance $d$ on $\dr C(M)$. 
Notice that a priori it is a pseudo-distance. 

Given $p\in\dr C(M)$ we 
construct a map $\iota: U\rightarrow \dr C(M)$, where $U$ is some open set
of $\mathbb H^2$, $\iota$ is bi-Lipschitz and preserves the length of curves, 
the image of $\iota$ is a neighbourhood of $p$.

We take a small neighbourhood $W$ of $p$ in $M$ and fix once and for all an
isometric identification with $W$ with some convex subset of $ H^3_1$.
Suppose without loss of generality that $p\in\dr_+ C(M)$. $\dr_+ C(M)\cap W$
can be regarded as the germ of a pleated surface of $H^3_1$, more
precisely we claim that there exists a complete 
convex in the past pleated surface in
$\mathbb H^3_1$, say $\Delta$, such that $\Delta\cap W'= \dr_+ C(M)\cap W'$
where $W'$ is a compact neighbourhood of $p$.

The existence of the map $\iota$ follows from the claim, thanks to the general
results about pleated surfaces in $H^3_1$ proved in~\cite{be-bo}.
Let us prove the claim. 

If $p\notin L$ then we can take $W'$ such that $W'\cap\dr_+ C(M)$ is totally
geodesic so the claim follows.

Suppose $p\in L$. We can choose $W_1\subset W$ that is pre-compact and such that
the leaves of $L$ meeting $W_1$ are exactly the leaves intersecting a small path
transverse to the leaf through $p$.

Consider the family $\mathcal F$ of spacelike planes of
$\mathbb H^3_1$ that are support planes of $\dr_+ C(M)\cap W$ at some point
$p\in W_1$.  The family $\mathcal F$ is pre-compact in $\mathbb H^3_1$ and  
there is a plane $P_0$ that does not intersect
any element of $\mathcal F$: indeed for a fixed $p_0\in\dr_+ C(M)$,
points of $P\in\mathcal F$ are connected to $p_0$ along spacelike geodesics so
it is sufficient to set $P_0$ to be the set of
points at distance $\pi/2$ from $p_0$. 

Now for $P\in\mathcal F$ denote by $C(P)$ the
convex set of $\mathbb H^3_1$ bounded by $P$ and $P_0$ and containing $\dr_+
C(M)\cap W$. Let $\Omega = \cap_{P\in\mathcal F} C(P)$. 
Let us enumerate some easy properties of $\Omega$.
\begin{enumerate}
\item $\Omega$
is a convex set of $\mathbb H^3_1$
\item $\Omega$ has two boundaries component. One of them is $P_0$. Let us set
  $\Delta$ to be the other component.
\item $\Delta\cap W_1=\dr_+ C(M)\cap W_1$.
\end{enumerate}
The last property is a consequence of the compactness of $\mathcal F$.

The last point to check is that we can choose  $W_1$ so that  
$\Delta$ is pleated. 

Suppose that for some $W_1$ 
some vertex occurs. By property (4) it is not difficult to see that there
are two bending lines in $\dr_+ C(M)$, say $l_1,l_2$ such that $l_i\cap
W\neq\varnothing$ and the geodesics of $\mathbb H^3_1$ extending $l_i\cap W$
meet each other at a point $q$. Consider points $p_i\in W\cap\ l_i$ and let
$T$ the geodesic triangle of $\mathbb H^3_1$ with vertices at $p_1,p_2$ and $q$.
Denote by $\hat l_i$ the segments joining $p_i$ to $q$.

Clearly $T\cap W_1$ embeds in $M$. 
Moreover the embedding $\sigma:T\cap W\rightarrow M$ extends on $(T\cap W)\cup
\hat l_1$ sending $\hat l_1$ on $l_1$ (and also to $(T\cap W)\cup \hat l_2$,
but a priori not on $(T\cap W)\cup\hat l_1\cup\hat l_2$).
If this isometric embedding
extends to an embedding on the whole $T$ we find a contradiction: the
image of $\hat l_i$ would be contained in $l_i$ and so the image of $q$ would
be contained in $l_1\cap l_2$.

Let $v_0$ the tangent vector at $p_1$ to the segment $p_1p_2$ and let
$v_t$ be the parallel transport of $v_0$ along $\hat l_1$. Consider the foliation 
of $T$ by geodesics arcs starting from $\hat l_1$ with direction $v_t$
(this is a foliation since $T$ is hyperbolic). Denote by $a_t$ the length of the segment $c_t$.

Let $v^*_t$ the parallel transport of $v_0$ along $l_1$ in $M$. Since the triangle 
$T$ does not embeds in $M$ there exists $t_0$ such that 
if $p_{t_0}$ is the corresponding point in $l_1$ then $\exp_{p_{t_0}} sv^*_{t_0}$ is 
defined for $s<b<a_{t_0}$, that means that $\exp_{p_{t_0}} b v^*_{t_0}$ is a singular point.

Notice that for a suitable choice of $W_1$ the factor $a_{t}$ can be close to $0$. 
On the other hand the vector $v_0$ runs in a compact set of $TM$ (independent of 
$W_1$), and so does the family $\{v^*_t|t>0\}$. Thus for any choice of  $W_1$
the pleated surface $\Delta$ contains some vertices there should be  a sequence 
in $L$ converging to some singular point. But
this implies that the closure of $L$ contains points on the singular locus and
this contradicts Lemma~\ref{lm:bd}.

Eventually we can choose $W'$ such that $\Delta$ is a complete pleated surface in
$\mathbb H^3_1$. Thus there is an isometry (that is a bijective map preserving
the distance) $B:\mathbb H^2\rightarrow\Delta$. Then $B: B^{-1}(\dr_+ C(M)\cap
W)\rightarrow\dr_+ C(M)\cap W$ is the isometry we are looking for.

In fact in~\cite{be-bo} 
the map $B$ is described in some more explicit way.
It is shown that there is a measured geodesic lamination $\lambda_\Delta$ on
$\mathbb H^2$ such that
\begin{enumerate}
 \item the bending locus of $\Delta$ is the image of the support of
 $\lambda_\Delta$;
 \item Up to post-composing with an isometry of $\mathbb H^3_1$ we have
\[
      B(x)=(\beta^R(x_0,x),\beta^L(x_0,x))I(x)
\]
where $x_0$ is a point, $\beta^R$ and $\beta^L$ are the right and left cocycles 
associated to $\lambda_\Delta$ as in Section~\ref{sec:earth} and $I$ is the standard embedding of 
 $\mathbb H^2$ in $H^3_1$ defined in Section \ref{sec:back}.
\item The lamination $\lambda_\Delta$ is determined by the bending: the
  bending of $\mathbb H^2$ along $\lambda$ and $\lambda'$ coincide on a
  neighbourhood $U$ iff $\lambda|_U=\lambda'|_U$.
\end{enumerate}
    
This show that it is possible to equip $L$ on $W\cap\dr_+ C(M)$ with a
transverse measure that is the image of the transverse measure on the
corresponding neighbourhood of $\mathbb H^2$.
Notice that by property (3) the transverse measures defined on different
neighbourhoods match on the intersecton, giving rise to a transverse measure
on $L$ on the whole of $\dr_+ C(M)$. Let $\lambda_+$ be the corresponding
lamination. Let $F_+$ be the hyperbolic structure on $\dr_+ C(M)$ and
$\lambda_+$ be the bending lamination. By point (3) these data determines the
developing map of $\dr_+ C(M)$ and thus the germ of the structure around
$\dr_+ C(M)$. From the uniqueness of the maximal extension, they determines the
whole of $M$.
\end{proof}

\begin{remark}
Note that the arguments given here are almost the same as for the
corresponding hyperbolic setting, as in the appendix of \cite{qfmp}.
Note also that the condition that the cone angles at the singularities are less
than $\pi$ seems to be really necessary to insure that the boundary of the
convex core is orthogonal to the singular lines. An interesting example 
can be found in \cite{benedetti-guadagnini}, it has ``particles'' with cone
angles equal to $\pi$ and it seems that the boundary of the convex core
is not ``orthogonal'' to those particles, and that it is bent along a 
geodesic segment joining its intersections with the two singular lines.
\end{remark}

\subsection*{Reconstruction from the boundary of the convex core.}

Thanks to Proposition~\ref{adstohyp:prop} there is a well-defined injective map
\[
   \cGH_{\Sigma,
   r,\theta}\rightarrow\cT_{\Sigma,r,\theta}\times\cML_{\Sigma,r}
\]
associating to $M$ the hyperbolic metric on the future boundary of the convex
core, say $h_+$, and the bending lamination, say $\lambda_+$.
The aim of this subsection is to show that this map is bijective, giving a
parametrization of $\cGH_{\Sigma,r,\theta}$ in terms of the embeding data of
the future boundary of the convex core (an analog parametrization is 
possible in terms of the embedding data of the past boundary of the convex core).

\begin{lemma}\label{lm:no}
Let $S$ be a pleated surface in $M$ convex in the past. Then there is no point
in the past of $S$ at distance $\pi/2$ from $S$.
\end{lemma}
\begin{proof}
Suppose by contradiction that there exists a point $p$ at distance
$\pi/2$. Let $r(p)\in S$ the point realizing the distance and let $P$ be the
set of points in $M$ that can be joined to $p$ by a timelike segment of length $\pi/2$.
$P$ is an immersed geodesic plane and clearly it is a support plane for $S$.
Moreover $S\cap P$ is convex in $P$ and without vertices. 
Thus the interior of $P\cap S$ is
contained in $S\setminus L$ whereas the boundary of $P\cap S$ is contained in
$L$. In particular $S\cap P$ contains a leaf, say $l$, of $L$.

If $l$ were closed, $l$ would be homotopic to the constant loop $p$ in $M$. 
Since $l$ is not trivial in $S$, it is not trivial in $M$ and this gives a
contradiction.

Suppose now that $l$ is open. We know there is a another leaf, $l'$, in
the closure of $l$ (the proof of this point can be done as for non-singular
hyperbolic surfaces, see e.g. \cite{FLP}). 
Moreover we can choose $l'$ that is not a boundary leaf --- that means that
support planes for points in $l'$ intersects $l'$ just in a geodesic segment.

Take a sequence $q_n\in l$ converging to 
$q_\infty\in l'$.  Timelike geodesics
connecting $p$ to $q_n$ go to a  geodesic $c$ connecting $p$
to $q_\infty$ that is not spacelike. If this geodesic segment is timelike 
then $q_\infty$ lies on $P$, that is the plane
orthogonal to the segment joining $p$ to $q_\infty$ is a support plane for
$q_\infty$ that contains many $q_n$. This contradicts the choice of $l'$.

Thus $c$ is lightlike. But this holds for every points of $l'$. On the other 
hand it is not difficult to prove that on a spacelike geodesic there are only 
a finite number of points 
connected to $p$ along a lightlike geodesic. This leads to a contradiction.
\end{proof}

\begin{prop} \label{lm:only}
Let $M$ be a GHM AdS manifold with particles. The only convex pleated surfaces
in $M$ are the future and past boundary components of the convex core. 
\end{prop}

\begin{proof}
If $S$ is a pleated surface convex in the past, it is contained in the future
of the convex core $C(M)$. Take $p\in S$ that is not in $C(M)$ and take a
point $q\in \dr_+C(M)$ such that $p\in I^+(q)$. Take a smooth convex surface
$S'$ in a neighbourhood of $\dr_+ C(M)$ and $q'\in S'$ such that $p\in
I^+(q')$. By Lemma~\ref{lm:existence}, $M$ contains a  timelike geodesic of length equal to
$\pi/2$ arriving at $q$. So there is a timelike path of length bigger than
$\pi/2$ arriving at $p$. But this contradicts Lemma~\ref{lm:no}
\end{proof}

\begin{prop} \label{lm:reconstruct}
Let $h$ be a hyperbolic metric with cone singularities (of angles $\theta_1,
\cdots, \theta_n\in (0,\pi)$) on $S$, and let $\lambda$ be a measured
bending lamination in the complement of the cone points. There is a unique 
GHM AdS metric with particles on $S\times (0,1)$ such that $h$ and
$\lambda$ are the induced metric and measured bending lamination on the
future boundary of the convex core.   
\end{prop}

\begin{proof}
The hyperbolic metric $h$ and the measured lamination $\lambda$ determine an
isometric embedding of the universal cover of the complement of
the cone points in $S$ into $AdS_3$ which is equivariant under an action of
the fundamental group $\Gamma$ of the complement of the cone points in $S$.
More precisely, the developing map can be explicitly written in terms of the
left and right cocycles  $\beta^l$ and $\beta^r$ associated to $\lambda$.
In fact 
\begin{equation}\label{eq:bend}
   dev(x)=(\beta^r(x_0,x),\beta^l(x_0,x)) I(dev_0(x))
\end{equation}
where $dev_0$ is the developing map of $h$ and $I$ is the standard embedding
$\mathbb H^2$ into $H^3_1$. 

The fact that $dev$ is locally injective and locally convex can be proved 
as in~\cite{benedetti-bonsante} in the non-singular case.
The only point to check is that it induces a hyperbolic structure on the surface $S$ with
cone singularities. On the other hand, since the singular locus is far from the 
lamination, the cocycles $\beta^r$ and $\beta^l$ are
trivial on $\pi^{-1}(U)$ for some neighbourhood  $U$ of the puncture. 
Thus on $\pi^{-1}(U)$ the map $dev$ is conjugated with $dev_0$, 
so the induced metric on $S$ looks like $h$ in a neighbourhood of a cone point.

Consider the normal exponential map
of $S$ towards the convex side of $S$. It is the map:
$$ G:N^1S\times (0,\pi/2)\rightarrow M~. $$
Here $N^1S$ is the unit normal bundle of $S$, i.e., the set of
unit vectors at points of $S$ for which the oriented orthogonal plane
is a support plane of $S$. The map is defined by sending $(n,t)$, where $n$
is a unit vector at $x\in S$, to $\exp_x(tn)$, where $\exp_x$ is the
exponential map at $x$.

The convexity of $S$ then shows that this map is locally injective on
$S\times (0,\pi/2)$ (as seen in the proof of Lemma \ref{dist:lem}). 
So this map can be used to pull back the AdS metric to a
locally AdS metric on $S\times (0,\pi/2)$, with cone singularities at the
lines $x\times (0,\pi/2)$, where $x$ is a cone point of $S$. This shows that 
$S$ has an embedding into an AdS manifold $N$ with image a convex pleated
surface with induced metric $h$ and measured pleating lamination $\lambda$.
By definition $N$ is contained in a GHMC AdS manifold $M$, also containing
a pleated surface with induced metric $h$ and measured 
pleating lamination $\lambda$.
 
In addition, Lemma \ref{lm:only} shows that $S$ can only be a connected
component of the convex core of $M$. Since a GHMC AdS manifold is obviously
determined by the future boundary of its convex core, the lemma follows.
\end{proof}

\subsection*{From the convex core to earthquakes.}

There is an important relation between convex pleated surfaces in GHMC AdS
manifolds and earthquakes on hyperbolic surfaces, which was discovered by Mess
\cite{mess}. It can be stated for convex cores of GHMC AdS manifolds with
particles as follows. 

\begin{lemma} \label{lm:diagram}
Let $M$ be a GHM AdS manifold with particles. Let $h_+, h_-$
be the induced metrics on the upper and lower boundaries of the
convex core, let $\lambda_+, \lambda_-$ be the measured bending
lamination of those upper and lower boundary components, and let
$\mu_l, \mu_r$ be the right and left hyperbolic metrics. Then
$$ \mu_l = E_l(\lambda_+)(h_+) = E_r(\lambda_-)(h_-)~, ~~
\mu_r = E_l(\lambda_-)(h_-) = E_r(\lambda_+)(h_+)~. $$
It follows that $\mu_l = E_l(2\lambda_+)(\mu_r) = E_r(2\lambda_-)(\mu_r)$. 
\end{lemma}

\begin{diagram}
& & h_+ & & \\
& \ldTo^{E_r
(\lambda_+)} & & \rdTo^{E_l(\lambda_+)} & \\
\mu_l & & & & \mu_r \\
& \luTo_{E_l(\lambda_-)} & & \ruTo_{E_r(\lambda_-)} & \\
& & h_- & & 
\end{diagram}

\begin{proof}
As stated in~\cite{minsurf}, the holonomy of $\mu_l$ is the projection 
of the holonomy of $M$ on the first factor of 
$Isom(\mathbb H^3_1)=PSL(2,\mathbb R)\times PSL(2,\mathbb R)$.
On the other hand, by formula~(\ref{eq:bend}), such a representation is simply
\[
   \rho_l(\gamma)=\beta_r(x_0,\gamma x_0)\rho_+(\gamma)
 \]
 where $\rho_+$ is the holonomy for $h_+$. Thus by formula~(\ref{eq:bend}), 
$\mu_+$ and $E_r (h_+ )$ share the same holonomy.
It is well known that, for hyperbolic surfaces with cone singularities of
angle less than $\pi$, the holonomy determines the hyperbolic structure
(this can be prove by the same argument, based on pant decomposition, as
for non-singular hyperbolic surfaces).
Since we are assuming cone angles less than $\pi$, the holonomy 
determines the structure, and we conclude that $\mu_+$ is equal to $E_r(h_+)$.
\end{proof} 

We can conclude from Lemma~\ref{lm:diagram} that Theorem~\ref{tm:mess} and 
Theorem~\ref{tm:earthquake} are equivalent. In fact, the composition 
\[
\begin{CD}
 \cT_{\Sigma,n,\theta} \times\cML_{\Sigma,n} @>I_+>>
\cGH_{\Sigma,n,\theta}@>(\mu_l,\mu_r)>>\cT_{\Sigma,n,\theta}\times\cT_{\Sigma,n,\theta}
 \end{CD}
 \]
 is the map
 \[
   (h,\lambda)\mapsto (E^r_\lambda(F), E^l_\lambda(F))
 \]
 Since $E^l_\lambda=(E^r_\lambda)^{-1}$, it is easy to see that the $(\mu_l,\mu_r)$ is bijective if 
and only if so is the map
 \[
   (h,\lambda)\mapsto (h,E^r_\lambda(h)),
 \]
and that in turn is equivalent to require that for all $h$ 
the map $E^r_\cdot(h):\cML_{\Sigma,n}\rightarrow\cT_{\Sigma,n,\theta}$ 
is bijective.

\section{Local deformations}

This section is devoted to local (or infinitesimal) deformations of 
GHM AdS manifolds with particles.  

\subsection*{Convex space-like surfaces in GHMC manifolds.}

In this section we consider a closed, convex space-like 
surface $S\subset M$ which is
orthogonal to the singularities. We call $S_r$ the regular set of $S$ --- the
complement in $S$ of the set of singular points --- and
$\Gamma:=\pi_1(S_r)$. There is a natural morphism
$$ \phi_M:\Gamma\rightarrow SO_0(2,2) $$
obtained from the holonomy representation of $M$, because
$\pi_1(S)=\pi_1(M_r)$. 

Given $S$ it is also possible to define two hyperbolic metrics on 
it, called $\mu_l$ and $\mu_r$ in the introduction. Note that this
depends on the fact that $S$ is convex (it is actually sufficient
to suppose that the curvature of the induced metric on $S$ does not
vanish, but this happens to be true for all convex surfaces). 

Those two metrics have cone singularities at the intersections of
$S$ with the singular lines of $M$, see \cite{minsurf}. Note that this point
depends on the fact that $S$ is orthogonal to the singular lines of $M$.

\subsection*{The left and right representations.}

We introduce a simple notation. For each cone point $x_i$ of $S$, $1\leq i\leq
n$, we call $\gamma_i$ the element of $\Gamma$ corresponding to the simple
closed curve going once around $x_i$.
 
A direct consequence of the Lemma \ref{lm:left-right} and of
the following remark is that the
images by $\phi_l$ and by $\phi_r$ of $\gamma_i$ 
are hyperbolic rotations of angle $\theta_i$, where $\theta_i$ is the 
angle of the singular curve of $M$ which intersects $S$ at $x_i$.

\begin{remark}\label{rk:rotation}
Let $\rho\in SO_0(2,2)$ be an AdS isometry, and let $\rho_l,\rho_r$ be its
left and right components. Let $\alpha\in (0,2\pi)$. $\rho$ is a (pure)
rotation of angle $\alpha$ around a time-like geodesic in AdS if and only if
$\rho_l$ and $\rho_r$ are both hyperbolic rotations of angle $\alpha$.  
\end{remark}

\subsection*{Deformations of the holonomy representation of $M$.}

\begin{lemma}\label{lm:representations}
The first-order deformations of $M$, among GHM AdS manifolds with particles 
with the same cone angles, are parametrized
by the first-order deformations of $\phi_M$, as a morphism from $\Gamma$ to
$SO_0(2,2)$, such that, for each $i\in \{ 1, \cdots, n\}$, the image of
$\gamma_i$ remains a pure rotation of angle $\theta_i$.
\end{lemma}

\begin{proof}
Consider a first-order deformation of the AdS metric on $M$, among GHMC AdS
metrics with the same cone angles. The corresponding first-order variation of
the holonomy representation of $M$ is then a first-order deformation of
$\phi_M$, and the images of the $\gamma_i$ remain a pure rotation of angle
$\theta$. 

Consider now a one-parameter deformation $\phi_t$ of $\phi_M$, $t\in
[0,\epsilon]$, such that the images of
each $\gamma_i$ remains a pure rotations of angle $\theta_i$. 
Let $\Mt$ be the universal cover of the complement of the singular lines in
$M$. There is natural local isometry from $\Mt$ to $AdS_3$, the developing map
of $M$.
Let $\St$ be the universal cover of the regular part of $S$. 
Then the developing map of $M$ restricts to an immersion:
$$ \psi_M: \St\rightarrow AdS_3~, $$
which is equivariant under the action $\phi_M$. Its image is a locally
convex surface, which is ramified at the images of the cone points.

Given the one-parameter deformation $\phi_t$, it is possible to construct
a one-parameter deformation $\psi_t$ of $\psi_M$, among  
embeddings of $\St$ into $AdS_3$, such that $\psi_t$ is equivariant 
under the action of $\phi_t$. This can be achieved for instance by
choosing a fundamental domain $D$ in $\St$ and constructing a deformation
of $\psi_M$ on $D$ in such a way that $D$ can be ``glued'' to its images
under the action of a set of genetators of $\Gamma$.
Moreover it is possible to choose this deformation so that, for $t$
small enough, the image of $\St$ remains locally convex. 

Then for $\alpha>0$ small enough and for $t$ small enough, we can consider the
normal exponential map:
$$ \exp_n:\St\times (-\alpha,\alpha)\rightarrow AdS_3~, $$
sending $(x,t)$ to the image of $tn_x$ by the exponential map at 
$\psi_t(x)$, where $n_x$ is the future-oriented unit normal orthogonal
to $\psi_t(\St)$ at $\psi_t(x)$. This map $\exp_n$ is a local homeomorphism,
so that it can be used to pull back the AdS metric of the target space to
an AdS metric on $\St\times (-\alpha,\alpha)$, which has a natural
isometric action of $\Gamma$ through $\phi_t$. 

The quotient $\St\times (-\alpha,\alpha)/\phi_t(\Gamma)$ is an AdS manifold
with particles, which contains a closed, locally convex, space-like 
surface (the quotient of $\psi_t(\St)$. So its maximal extension is a GHMC AdS
manifold with particles, with holonomy representation equal to $\phi_t$, as
needed. 

Note that any one-parameter deformation of $M$ (still under the same angle
conditions) can be constructed in this manner, and that the resulting
manifold depends only on the variation of the holonomy representation.
This completes the proof of the lemma.
\end{proof}

\subsection*{A key infinitesimal rigidity lemma.}

We now have the tools necessary to prove the main lemma of this section: the
first-order deformations of $M$ are parametrized by the first-order
deformations of its left and right hyperbolic metrics. 

\begin{lemma}\label{lm:localrig}
The map $(\mu_l,\mu_r):\cGH_{\Sigma,n,\theta}\rightarrow
\cT_{\Sigma,n,\theta}\times \cT_{\Sigma,n,\theta}$ is a local
homeomorphism.
\end{lemma}

\begin{proof}
According to Lemma \ref{lm:representations}, the first-order deformations of
$M$ are parametrized by the first-order deformations of its holonomy
representation, among morphisms of $\Gamma$ in $SO_0(2,2)$ sending each
$\gamma_i$ to a pure rotation of angle $\theta_i$. But Remark
\ref{rk:rotation} shows that those deformations are characterized by the 
deformations of the left and right hyperbolic metrics, 
which can be any deformations
of $\mu_l$ and $\mu_r$ among hyperbolic metrics with the same angle at the
cone singularities of $S$. 
\end{proof}

\subsection*{Consequences for earthquakes.}

The previous lemma has a direct application for earthquakes on hyperbolic
surfaces with cone singularities. 

\begin{lemma}  \label{lm:localrig-2}
Let $h_0\in \cT_{\Sigma,n,\theta}$. The map $E^r_\cdot(h_0):\lambda\mapsto
E^r_\lambda(h)$ is a local homeomorphism.
\end{lemma}

\begin{proof}
Fix $\lambda_0$.
Let $I:\cML_{\Sigma,n}\times\cT_{\Sigma,n,\theta}\rightarrow\cGH_{\Sigma,n,\theta}$ be the parameterization given by 
Lemma \ref{lm:diagram}. Since $I$ is continuous, thanks to 
Lemmas \ref{lm:localrig} and \ref{lm:diagram}, we can take a neighbourhood $U$ of $(h_0,\lambda_0)$
such that the map
\[
   (h,\lambda)\mapsto (E^r_{\lambda/2}(h), E^l_{\lambda/2}(h))
 \]
 is injective on $U$.
 The set  $V$  of laminations $\lambda$ such that $(h_0,\lambda)\in U$
 is a neighbourhood of $(h_0,\lambda_0)$ and both $E^r_\cdot(h_0)$ and $E^l_\cdot(h_0)$ are injective on $V$.
\end{proof}
 
\section{Compactness}

This section is devoted to the proof of Lemma \ref{lm:compact}, which states
that, for a fixed element $\mu\in \cT_{\Sigma,n,\theta}$, the map
$E^r_\cdot(\mu):\cML_{\Sigma,n}\rightarrow 
\cT_{\Sigma,n,\theta}$ is proper. In the whole section we fix
$\theta=(\theta_1, \cdots, \theta_n)\in [0,\pi)^n$.

\subsection*{A compactness lemma for earthquakes.}

\begin{lemma}\label{lm:compact:est}
Given $\lambda\in\cML_{\Sigma,n}$ and  $\mu\in\cT_{\Sigma,n,\theta}$ 
let $\mu'=E^r_{\lambda}(\mu)$. Then, for every closed geodesic $\gamma$ of
$\Sigma$ the following estimate holds
\[
    \ell_\mu(\gamma)+\ell_{\mu'}(\gamma)\geq\lambda(\gamma)
\]
where $\ell_\mu(\gamma)$ denotes the length of $\gamma$ with respect to $\mu$.
\end{lemma}

\begin{proof}
By a standard approximation argument, it is sufficient to prove the statement
under the hypothesis that  $\lambda$ is   a weighted multicurve. Moreover we
can assume $\lambda(\gamma)>0$.

Let $\tilde\Sigma_\mu$ and $\tilde\Sigma_{\mu'}$ the metric universal covering
of the regular part of $\Sigma$ with respect 
$\mu$ and $\mu'$ respectively (here $\Sigma$ is
regarded as a punctured surface). The lamination $\lambda$ lifts to a
lamination $\tilde\lambda$ of $\tilde\Sigma_\mu$ and the right earthquake
along $\tilde\lambda$, say $\tilde E^r$, is the lifting of $E^r$.
Thus for every covering transformation $g$ of $\tilde\Sigma_\mu$ there exists
a unique covering transformation, say $H(g)$, such that the following
equivariance formula holds
\[
        \tilde E^r\circ g=H(g)\circ\tilde E^r\,.
\]

Let $g$ be a covering transformation of $\tilde\Sigma_\mu$ representing a loop
of $\Sigma$ freely homotopic to $\gamma$. There exists a 
$g$-invariant complete geodesic $A=A(g)$ in $\tilde\Sigma_\mu$ such
that the projection of $A(g)$ on $\Sigma$ is the
geodesic representative of $\gamma$ with respect to $\mu$.

Analogously the projection of $A(H(g))\subset\tilde\Sigma_{\mu'}$, is the 
geodesic representative of $\gamma$ with respect to $\mu'$.
The inverse image, $A'$, through $\tilde E^r$ of $A(H(g))$ is a $g$-invariant union
of disjoint geodesic segments whose end-points lie on some leaves of
$\tilde\lambda$. 
More precisely if $\{l_i\}_{i\in\mathbb Z}$ is the set of geodesics cutting $A$
enumerated so that $l_i$ and $l_{i+1}$ intersect $A$ in consecutive points
$p_i$ and $p_{i+1}$ of $A\cap\tilde\lambda$, then $A'$ is the union of 
geodesic segments joining a
point $q_i\in l_i$ to a point $r_{i+1}\in l_{i+1}$.

Let $A$ be oriented  in such a way that $g$ is a positive translation. 
The sequence $p_i$ can be supposed to be increasing. 
Moreover each $l_i$ can be oriented in
such a way that the intersection of $A$ with it is positive.
Let $x_i$ (resp. $y_i$) denote the signed distance of $q_i$ (resp. $r_i$)
from $p_i$ on $l_i$. Since after the right earthquake $A'$ becomes a
continuous line, $x_i-y_i$ is equal to the weight of the leaf $l_i$.

\begin{figure}
\begin{center}
\input geo.pstex_t
\end{center}
\end{figure}

If $g(p_0)=p_n$ for some $n>0$ clearly we have $x_{i+n}=x_i$ and
$y_{i+n}=y_i$. Moreover  $\ell_\mu(\gamma)$ is equal to the sum of the 
lengths of the geodesic segments $[p_0,p_1],\ldots, [p_{n-1},p_n]$, whereas
$\ell_{\mu'}(\gamma)$ is equal to the sum of the geodesic segments 
$[q_0,r_0],\ldots,[q_{n-1},r_n]$. From the triangular inequality we get that
\[
    x_i\leq y_{i+1}+\ell([p_i,p_{i+1}])+\ell([q_{i}, r_{i+1}])
\]
that is
\[
   x_i-y_{i+1}\leq \ell([p_{i},p_{i+1}])+\ell([q_i, r_{i+1}])\,.
\]
Summing the last inequality for $i=0,\ldots, n-1$ we get
\[
  \sum_{i=0}^{n-1} x_i-y_i\leq \ell_\mu(\gamma)+\ell_{\mu'}(\gamma)
\]
Since the left hand of this inequality is the mass of $\gamma$ with respect to
$\lambda$, the proof is complete.
\end{proof}

\subsection*{Proof of Lemma \ref{lm:compact}.}
Let $(\lambda_k)_{k\in\mathbb N}$ be a divergent sequence in 
$\cML_{\Sigma,n}$ and 
$\mu_k=E^r_{\lambda_k}(\mu)$ for some fixed $\mu\in\cT_{\Sigma,n,\theta}$.
We have to prove that $(\mu_k)_{k\in\mathbb N}$ is a divergent sequence in
$\cT_{\Sigma,n,\theta}$.

Since $(\lambda_k)_{k\in\mathbb N}$ is divergent, there exists a closed
geodesic $\gamma$ such that 
$\lambda_n(\gamma)\rightarrow+\infty$. 
Then by Lemma~\ref{lm:compact:est}, $\ell_{\mu_k}(\gamma)\rightarrow+\infty$,
so $(\mu_k)_{k\in\mathbb N}$ does not admit a convergent subsequence.

\section{Proof of the main results}

\subsection*{Proof of Theorem \ref{tm:earthquake}.}

As mentioned in section 2 we fix $\theta=(\theta_1, \cdots,
\theta_n)\in (0,\pi)^n$ and $h\in \cT_{\Sigma,n,\theta}$. We then consider 
the map $E^r_\cdot(h):\cML_{\Sigma,n}\rightarrow \cT_{\Sigma,n,\theta}$. It is
a local 
homeomorphism by Lemma \ref{lm:localrig-2}, and is proper by Lemma
\ref{lm:compact}. Therefore it is a covering. However $\cT_{\Sigma,n,\theta}$
is 
simply connected and $\cML_{\Sigma,n}$ is connected, therefore this map is a
homeomorphism. 

\subsection*{Proof of Theorem \ref{tm:mess}.}

Again we consider a fixed choice of $\theta=(\theta_1, \cdots,
\theta_n)\in (0,\pi)^n$. 
Let $\mu_l,\mu_r\in \cT_{\Sigma,n,\theta}$. By Theorem \ref{tm:earthquake}
there exists a unique $\lambda\in \cML_{\Sigma,n}$ such that
$\mu_r=E^r_\lambda(\mu_l)$. 

Let $\lambda_+:=\lambda/2$, and let $h_+:=E^r_{\lambda_+}(\mu_l)$. 
By Lemma \ref{lm:reconstruct} there exists a (unique) GHMC AdS metric
$g$ on $\Sigma \times (0,1)$ 
for which the induced metric and the measured bending lamination on the upper
boundary of the convex core are $h_+$ and $\lambda_+$, respectively.
It then follows from Lemma \ref{lm:diagram} that the left and right
hyperbolic metrics of $g$ are $\mu_l$ and $\mu_r$,
respectively. 

Conversely, let $g'$ be a GHMC AdS metric on $\Sigma \times (0,1)$ 
for which the left and right hyperbolic metrics are $\mu_l$ and $\mu_r$,
respectively. Let $h'_+$ and $\lambda'_+$ be the induced metric and 
measured bending lamination on the upper component of the boundar of
the convex core of $(\Sigma\times (0,1),g')$. Then Lemma \ref{lm:diagram}
shows that $\mu_r=E^r_{\lambda_+}(h'_+)$, while $\mu_l=E^l_{\lambda_+}(h'_+)$,
so that $\mu_r=E^r_{2\lambda_+}(\mu_l)$. It follows that $g'$ is the 
metric $g$ constructed above, and this finishes the proof of Theorem
\ref{tm:mess}.

\section{Some concluding remarks}

\subsection*{Reconstructing a GHMC AdS manifold from its convex core.}

The arguments developed above show that, given the convex core of a GHMC AdS
manifold $M$, it is possible to understand the global geometry of $M$
is a simple way. This is an immediate extension of statements already
well-known in the non-singular case, see \cite{benedetti-bonsante}.

\begin{lemma}
Let $M$ be a GHM AdS manifold with particles, which is topologically
$\Sigma\times (0,1)$, with cone angles $\theta_1, \cdots, \theta_n\in
(0,\pi)$. Let $\Omega_+$ be the set of points at distance at most $\pi/2$
in the past of the future boundary of the convex core, and let 
$\Omega_-$ be the set of points at distance at most $\pi/2$ in 
the future of the past boundary of the convex core. Then
$$ M = \Omega_+\cup \Omega_-~,~~ \Omega_+\cap\Omega_- = C(M)~. $$
Moreover,
$$ \vol(M) + \vol(C(M)) = \frac{\pi}{2}\left(2\pi\chi(\Sigma) + \sum_{i=1}^n
  (2\pi-\theta_i)\right) + \frac{L(\lambda)}{2}~, $$
where $\lambda$ is the measured bending lamination of the boundary of
the convex core and $L(\lambda)$ is its length.
\end{lemma}

Note that the quantity $2\pi\chi(\Sigma) + \sum_{i=1}^n (2\pi-\theta_i)$ is
$2\pi$ times a natural ``Euler characteristic'' of a closed surface with
cone singularities, it is equal for instance to the area of any hyperbolic
metric with prescribed singular angles on such a surface.

\begin{proof}
Let $\dr_-C(M)$ and $\dr_+C(M)$ be the past and future boundary components of
$C(M)$, respectively. Since $\dr_-C(M)$ is a locally convex surface (with 
the convex part of its complement in the future direction) we can consider
the normal exponential map of $\dr_-C(M)$, as in the proof of the previous
lemma. Again it is locally injective on time $t\in (0,\pi/2)$, and can be used
to pull back the AdS metric to a locally convex AdS metric (with particles) on
a ``slice'' of width $\pi/2$ in the future of $\dr_-C(M)$.

This construction, and the definition of $M$ as a {\it maximal}
globally hyperbolic space, shows that $M$ contains the image of this map,
which corresponds to the space $\Omega_+$ appearing in the lemma. 

On the other hand, there is no point in $M$ which is at distance larger than
$\pi/2$ in the future of $\dr_-C(M)$. Indeed, suppose that some point $x\in M$
is at distance $\pi/2$ in the future of $\dr_-C(M)$. Let $\gamma_0$ be a
maximizing geodesic in $M$ from $x$ to a point $y\in \dr_-C(M)$. $y$ is
contained in a maximal totally geodesic stratum $\sigma$ of $\dr_-C(M)$, 
and, for all points $y'\in \sigma$, the geodesic orthogonal to $\sigma$ at
$y'$ arrives at $x$ after time exactly $\pi/2$. 

But the universal cover $\tilde{\sigma}$
of $\sigma$ is non-compact, and we can consider a
sequence $(y_n)$ of points in it, going to infinity. Let $\gamma_n$ be the
projection to $M$ of the 
geodesic segment or length $\pi/2$ orthogonal to $\tilde{\sigma}$ at $y_n$, so 
that the other endpoint of $\gamma_n$ is $x$. Then the sequence $(\gamma_n)$
converges, in all compact sets containing $x$, to a light-like geodesic which
does not intersect the universal cover of $\dr_-C(M)$, contradicting the
global hyperbolicity of $M$. 

This proves that the future of $\dr_-C(M)$ is equal to the image by $G$ of 
$\dr_-C(M)\times (0,\pi/2)$, or in other terms $\Omega_+$. The same
argument applies for the past of $\dr_+C(M)$, and this proves 
the first part of the lemma.

For the second point, 
let $\lambda_+$ and $\lambda_-$ be the measured bending lamination of the
future and past boundary components of $C(M)$. The statement will clearly
follow if we prove that 
$$ \vol(\Omega_+) = \frac{\pi}{4}A(\dr_-C(M)) + \frac{L(\lambda_-)}{2}~, $$
where $A(\dr_-C(M))$ is the area of the induced metric on $\dr_-C(M)$ and 
$L(\lambda_-)$ is the length of the measured lamination $\lambda_-$. 
We will prove that this relation holds when the support of $\lambda_-$ 
is a disjoint union of closed curves; the general result for $\Omega_+$ then
follows by approximating $\dr_-\Omega_+$ by a sequence of pleated surfaces
with such a measured bending lamination.

So suppose that the support of $\lambda_-$ is the union of closed curves
$\gamma_1, \cdots, \gamma_p$, each with a weight $\lambda_i, 1\leq i\leq p$. 
Then $\Omega_+$ can be decomposed as the union of two areas:
\begin{itemize}
\item the set $\Omega_0$ 
of points $x\in \Omega_+$ which ``project'' to a point of
  $\dr_-C(M)$ which is not in the support of $\lambda_-$,
\item the sets $\Omega_i$ of points $x\in \Omega_+$ which project to
  $\gamma_i$. 
\end{itemize}
The volume of the first domain can be computed by integrating over the
distance to $\dr_-C(M)$, it is equal to:
$$ \vol(\Omega_0) = \int_{\dr_-C(M)\setminus \cup_i \gamma_i}
\int_{r=0}^{\pi/2} \cos(r)^2 dr da = \frac{\pi}{4} 
\int_{\dr_-C(M)\setminus \cup_i \gamma_i}
da = \frac{\pi}{4}A(\dr_-C(M))~. $$
The same kind of computation shows that: 
$$ \vol(\Omega_i) = L(\gamma_i) \lambda_i \int_0^{\pi/2}\cos(r)\sin(r) dr = 
\frac{L(\gamma_i) \lambda_i}{2}~, $$
and it follows that
$$ \vol(\Omega_+) = \frac{\pi}{4}A(\dr_-C(M)) + \frac{L(\lambda_-)}{2}~, $$
as needed.
\end{proof}

\subsection*{A note on the definition of GHMC manifolds used here.}

The reader might wonder why we consider here {\it convex} globally hyperbolic
manifolds, i.e., AdS manifolds which contain a space-like surface which is
convex (see Definition \ref{df:convex}). It is quite possible that a weaker
condition -- the existence of a compact Cauchy surface, which is not necessarily
convex -- is sufficient, and that any AdS manifold containing such a Cauchy
surface contains one which is convex. We do not further consider this question
here since it is quite distinct from our main centers of interest. 

\subsection*{Other possible proofs.}

It appears quite possible that arguments close to those used by Kerckhoff
\cite{kerckhoff} can be applied to the setting of hyperbolic surfaces with
cone singularities, however the extension is not completely clear. We suppose
that the condition that the cone angles are less than $\pi$ should appear also
in such arguments. 

Other arguments used without on non-singular surfaces, however, make a
stronger use of the geometry of the universal cover of the surface. Those are
presumably well adapted to the singular surfaces considered here. 

\bibliographystyle{amsplain}
\bibliography{../outils/biblio}

\end{document}